\documentclass[10pt]{amsart}

\usepackage{eucal}
\usepackage{amscd}
\usepackage[all,cmtip]{xy}
\usepackage{rotating}
\usepackage{hyperref,mathrsfs}
\usepackage{tikz}
\usepackage{amsmath,amssymb}

\newtheorem{theorem}{Theorem }[section]

\newtheorem{lemma}[theorem]{Lemma}

\newtheorem{corollary}[theorem]{Corollary}
\newtheorem{proposition}[theorem]{Proposition}
\newtheorem{question}[theorem]{Question}
\theoremstyle{definition}

\newtheorem{remark}[theorem]{Remark}
\newtheorem{observ}[theorem]{Observation}

\def\PG{\mathbf{PG}}

\def\eop{\hspace*{\fill}{\footnotesize$\blacksquare$}}
\newcommand{\Aut}{\mathrm{Aut}}
\newcommand{\id}{\mathrm{id}}

\newcommand{\wis}[1]{{\text{\em \usefont{OT1}{cmtt}{m}{n} #1}}}

\newcommand{\A}{\mathbb{A}}
\newcommand{\hP}{\mathbb{P}}

\newcommand{\Fun}{\mathbb{F}_1}
\newcommand{\Spec}{\wis{Spec}}
\newcommand{\Z}{\mathbb{Z}}
\newcommand{\Proj}{\wis{Proj}}

\newcommand{\bG}{\mathbb{G}}
\newcommand{\Hom}{\wis{Hom}}
\newcommand{\Adj}{\wis{Adj}}

\newcommand{\mF}{\mathcal{F}}

\newcommand{\fP}{\mathbf{P}}

\newcommand{\mX}{\mathcal{X}}

\newcommand{\F}{\mathbb{F}}

\newcommand{\mY}{\mathcal{Y}}

\newcommand{\PGL}{\mathbf{PGL}}
\newcommand{\PGaL}{\mathbf{P\Gamma L}}

\newcommand{\prf}{\textit{Proof. }}
\newcommand{\Sch}{\wis{Sch}}

\usetikzlibrary{matrix}
\usetikzlibrary{arrows}
\usepackage{pgf}
\usepackage{lscape}
\usepackage{listings}

\usepackage{enumitem}

\title[Automorphisms of Deitmar schemes]{Automorphisms of Deitmar schemes, I. Functoriality and trees}

\keywords{Field with one element; Deitmar scheme; loose graph; loose tree; zeta function; functoriality; automorphism group; dichotomy}

\author{Manuel M\'{e}rida-Angulo and Koen Thas}

\address{{Ghent University},
{Department of Mathematics},
{Krijgslaan 281, S25, B-9000 Ghent, Belgium,}
\texttt{manmerang@gmail.com}
\texttt{koen.thas@gmail.com}}

\date{}

\begin{document}
\maketitle

\begin{abstract}
In a recent paper \cite{MMKT}, the authors introduced a map $\mF$ which associates a Deitmar scheme (which is defined over the field with one element, denoted by $\Fun$) with any given graph $\Gamma$. By base extension, a scheme $\mX_k = \mF(\Gamma) \otimes_{\Fun} k$ over any field $k$ arises. In the present paper, we will show that all these mappings are functors, and we will use this fact to study automorphism groups of the schemes $\mX_k$. Several automorphism groups are considered: combinatorial, topological, and scheme-theoretic groups, and also groups induced by automorphisms of the ambient projective space.
When $\Gamma$ is a finite tree, we will give a precise description of the combinatorial and projective groups, amongst other results. 
\end{abstract}

\setcounter{tocdepth}{1}
\tableofcontents

\medskip
\section{Introduction}

In a recent paper \cite{MMKT}, the authors of the present text have introduced a new zeta function for finite graphs and a generalization of the latter, called
``loose graphs.'' Loose graphs look like graphs but it is allowed that edges have $1$ and even $0$ vertices. 
In fact, what the authors did was associate an extended Deitmar scheme $\mF(\Gamma)$ to each such loose graph $\Gamma$, and then show that 
the obtained Deitmar schemes enjoy a number of properties which allow us to attach a  {\em Kurokawa zeta function} to the scheme. Deitmar schemes 
are schemes defined over the field with one element, $\Fun$, and are the main objects in the algebraic geometry of monoids. Understanding such 
schemes was one of the main motivations to start our study. 

A second driving force after this paper is the following. According to the map $\mF$, to each vertex in a loose graph $\Gamma$ is associated 
an affine space of dimension the degree of the vertex, and the affine spaces as such obtained are glued according to rules which can be read from the 
incidences in $\Gamma$. A natural question then becomes what features of the schemes $\mF(\Gamma) \otimes_{\Fun} k$ can be determined directly 
from $\Gamma$ (that is, can be expressed in terms of degrees, cycles, etc.)?

\subsection*{Functoriality} 
For each field $k$, one obtains a $k$-scheme $\mX_k = \mF(\Gamma) \otimes_{\Fun} k$
and the schemes $\mX_k$ come with the same Kurokawa zeta function. Although it is mentioned in \cite{MMKT}, it was not shown that the association
\begin{equation}
\mF: \Gamma \mapsto \mF(\Gamma)
\end{equation}
defines a covariant functor from loose graphs to extended Deitmar schemes, and this is the first goal of the present paper: showing that for each field $k$, including $\Fun$, the map 
\begin{equation}
\mF_k: \Gamma \mapsto \mF(\Gamma)\otimes_{\Fun}k
\end{equation}
is a functor (where $\mF_{\Fun} = \mF$). (An easy but very interesting feature of this part of the paper is that morphisms which are not monomorphisms 
are in nature not of ``$\Fun$-type,'' because they tend to introduce additions. A detour to $\mathbb{F}_2$ is needed to solve this problem.)

\subsection*{Automorphisms}

Next, we study schemes coming from (loose) trees, and in particular, we determine the automorphism groups of these schemes (over any field) in terms of data associated to the (loose) trees. One of the main results is the following.

\begin{theorem}
\label{MTtreesintro}
Let $T$ be a loose tree, and let $k$ be any field. Put $\mX_k = \mF(T) \otimes_{\Fun} k$, and consider the embedding
\begin{equation}
\iota: T \ \hookrightarrow\ \mX_k.
\end{equation}
Let $I$ be the set of inner vertices of $T$, and let $T(I)$ be the subtree of $T$ induced on $I$. 
We have $\PGaL(\mX_k) = \Aut^{\mathrm{proj}}(\mX_k)$ is isomorphic to 
\begin{equation}
\Big(\Big(\prod^{\mathrm{centr}}_{w \in I}S(w)\Big) \rtimes \Aut(T(I))\Big) \rtimes k^{\times}.
\end{equation}
\end{theorem}

Here, $\Aut^{\mathrm{proj}}(\mX_k)$ denotes the automorphism group of $\mX_k$ which is induced by the automorphism group of the ambient projective space.
(For other notational details we refer the reader to the body of the text.)

Three different automorphism groups are considered: the aforementioned ``projective group,'' the {\em combinatorial automorphism group}, which is the automorphism group of $\mX_k$ considered as an incidence geometry, and the {\em topological automorphism group}, which is the group of homeomorphisms  $\mX_k \mapsto \mX_k$ (where the topology is that coming from the scheme). 

Many other results are obtained in this context, and we study some particular examples in much detail before passing to general theorems.

\subsection*{Edge-relation dichotomy}

The reader must note that calculating automorphism groups of schemes related to (loose) graphs is a problem of high complexity, since one needs to calculate the automorphism group of the (loose) tree {\em before} being able to handle the associated schemes. Therefore, we introduce the ``inner graph property,'' which is a property for (loose) graphs which leads to more or less direct calculation of the automorphism groups of the schemes in function of the groups of the (loose) graphs.

Finally, we study the ``edge-relation dichotomy,'' a phenomenon which predicts {\em when} a loose graph {\em has} the inner graph property, in function of its distance to the ambient space and to a tree. 

\section{Deitmar schemes}

For the definition of Deitmar schemes one has to recall first some important definitions and constructions. We define an {\em $\Fun$-ring} $A$ to be a multiplicative commutative monoid with an extra absorbing element 0.
Let $\Spec(A)$ be the set of all {\em prime ideals} of $A$ together with a Zariski topology. We refer to \cite{Deitmarschemes2} for the definition of prime ideals of a monoid. This topological space endowed with a structure sheaf of $\Fun$-rings is called an {\em affine Deitmar scheme} in the same way as affine schemes defined over a field $k$, or $\mathbb{Z}$. We define a {\em monoidal space} to be a pair $(X, \mathcal{O}_X)$ where $X$ is a topological space and $\mathcal{O}_X$ is a sheaf of $\Fun$-rings defined over $X$. A {\em Deitmar scheme} is then a monoidal space such that for every point $x\in X$ there exists an open subset $U \subseteq X$ such that $(U, \mathcal{O}_X|_{U})$ is isomorphic to an affine Deitmar scheme.\medskip

\subsection{Affine space} One of the most important examples of Deitmar schemes is the affine space $\mathbb{A}^n_{\Fun}$. Let us describe its construction.\medskip

Define the {\em monoidal ring on $n$ variables} $X_1, \ldots, X_n$ as the monoid
\begin{equation}
\Fun[X_1, \ldots, X_n]:=\{0\}\cup \{X_1^{a_1}\cdots X_n^{a_n} ~|~ a_i \in \mathbb{N} \},
\end{equation}

\noindent i.e, the union of $\{0\}$ and all the monomials generated by the variables $X_i$. Let us call $A:=\Fun[X_1, \ldots, X_n]$; then the {\em n-dimensional affine space over $\Fun$} is defined as the monoidal space $\Spec(A)$ and denoted by $\mathbb{A}^n_{\Fun}$. Note that all the prime ideals of $A$ are finite unions of ideals of the form $(X_i)$, where $(X_i)=\{X_ia ~|~ a\in A \}$. For a more detailed definition of Deitmar schemes and the structure sheaf of $\Fun$-rings, we refer to \cite{Deitmarschemes2}, or \cite{Chap2}.\medskip

\subsection{Congruence schemes} A more general version of Deitmar scheme is a so-called {\em congruen\-ce scheme}. For the definition of congruen\-ce scheme, we refer to \cite{Deitmarcongruence}. Let us just mention that congruence schemes are defined in terms of {\em sesquiads}. A sesquiad is a monoid $A$ endowed with an {\em addition} or $+$-{\em structure}; this $+$-structure allows addition for a certain set of elements in the monoid $A$. It is known that the category of monoids is a full subcategory of the category of sesquiads.\medskip

A sesquiad is said to be {\em integral} if $1\neq 0$ and 

\begin{equation*}
af=bf \implies (a=b ~ \mbox{ or } ~ f=0).
\end{equation*}

A {\em congruence} on a sesquiad $A$ is an equivalence relation $\mathcal{C}\subseteq A\times A$ such that there is a sesquiad structure on $A/\mathcal{C}$ making the projection $A \rightarrow A/\mathcal{C}$ a morphism of sesquiads. If $A/\mathcal{C}$ is integral, the congruence $\mathcal{C}$ is called {\em prime}. We denote by $\Spec_{c}(A)$ the set of all prime congruences on the sesquiad $A$ with the topology generated by all sets of the form 
\begin{equation*}
D(a,b)=\{ \mathcal{C} \in \Spec_{c}(A) ~|~ (a,b)\notin \mathcal{C }\}, ~~~~~~~a,b \in A.
\end{equation*}

In a similar way as for monoids, one can define a {\em structure sheaf of sesquiads} and a {\em sesquiaded space}. We call an {\em affine congruence scheme} to be a sesquiaded space that is of the form $(\Spec_{c}(A), \mathcal{O}_A)$, for $A$ a sesquiad and $\mathcal{O}_A$ its corresponding structure sheaf, and a {\em congruence scheme} to be a sesquiaded space $X$ that locally looks like an affine one. 

\subsection{The $\Proj_{c}$-construction} 

Consider the monoid $\Fun[X_0, X_1, \ldots , X_m]$, where $m\in\mathbb{N}$ and see it as a sesquiad together with the {\em trivial addition}. Since any polynomial is homogeneous in this sesquiad, we have a natural grading
\begin{equation*}
\Fun[X_0, \ldots , X_m]= \bigoplus_{i\geq 0}R_i=\coprod_{i\geq 0}R_i,
\end{equation*}

\noindent where $R_i$ consists of the elements of $\Fun[X_0, \ldots, X_m]$ of total degree $i$, for $i\in \mathbb{N}$. We defined then the {\em irrelevant congruence} as 
\begin{equation*}
\mbox{Irr}_{c}= \langle X_0 \sim 0 , \ldots, X_m \sim 0 \rangle.
\end{equation*}

Now we can proceed with the usual $\Proj$-construction of projective schemes. We define $\Proj_{c}(\Fun[X_0,\ldots, X_m])$ as the set of prime congruences of the sesquiad $\Fun[X_0,\ldots,X_m]$ which do not contain $\mbox{Irr}_{c}$. The closed sets of the topology on this set are defined as usual: for any $(a,b)$ pair of elements of $\Fun[X_0,\ldots,X_m]$, we define
\begin{equation*}
V(a,b):=\{\mathcal{C}~|~ \mathcal{C}\in \Proj_{c}(\Fun[X_0,\ldots, X_m]) , ~~ a\sim_{\mathcal{C}} b\},
\end{equation*}
\noindent and these sets form a basis for the closed set topology. Defining the structure sheaf as in \cite{Deitmarcongruence}, one obtains that $\Proj_c(\Fun[X_0,\ldots,X_m])$ is a {\em projective} congruence scheme, and it is this scheme which will be used in this paper for a projective $\Fun$-space. Its closed points naturally correspond to the $\mathbb{F}_2$-rational points of the projective space $\mathbb{P}^{m}(\mathbb{F}_2)$ (but the latter has a finer subspace structure, and also a different algebraic structure).

\section{Loose graphs and Deitmar schemes}

In this section we will briefly describe how one can associate a Deitmar scheme to a {\em loose graph} $\Gamma$, which is a graph in which edges with 0 and 1 end points are also allowed, through a functor, which we call $\mathcal{F}$. This functor must obey a set of rules, namely: \medskip

\begin{itemize}
\item[COV]
If $\Gamma \subset \widetilde{\Gamma}$ is a strict inclusion of loose graphs, $\mF(\Gamma)$ also is a proper subscheme
of $\mF(\widetilde{\Gamma})$.
\item[LOC-DIM]
If $x$ is a vertex of degree $m \in \mathbb{N}^\times$ in $\Gamma$, then there is a neighborhood $\Omega$ of $x$ in 
$\mF(\Gamma)$ such that $\mF(\Gamma)_{\vert \Omega}$ is an affine space of dimension $m$.
\item[CO]
If $K_m$ is a sub complete graph on $m$ vertices in $\Gamma$, then $\mF(K_m)$ is a closed sub projective space 
of dimension $m - 1$ in $\mF(\Gamma)$. 
\item[MG]
An edge without vertices should correspond to a multiplicative group.
\end{itemize}
\medskip

Rule (MG) implies that we have to work with more general Deitmar schemes since the multiplicative group $\mathbb{G}_m$ over $\Fun$ is defined  to be isomorphic to
\begin{equation*} \Spec(\Fun[X, Y]/(XY=1)),\end{equation*} where the last equation generates a congruence on the free abelian monoid $\Fun[X_1, X_2]$. The reader can find a more detailed explanation of this association in \cite{MMKT}.

\section{Functoriality}

We will prove in this section that the above association, called $\mathcal{F}$, is indeed a proper functor. For this purpose, it remains to verify that morphisms of loose graphs induce morphisms of congruence schemes.\medskip

\subsection{Local action}Let $\Gamma$ be a finite loose graph and $f$ a loose graph automorphism of $\Gamma$. Remember that $\mathcal{F}(\Gamma)$ is the union of finite dimensional affine schemes defined from the vertices of $\Gamma$, i.e., 
\begin{equation*}
\displaystyle \mathcal{X} = \mathcal{F}(\Gamma) = \bigcup_{v\in V(\Gamma)} \Spec(A_v)
\end{equation*}

\noindent where $A_v$ is a finite $\Fun$-ring isomorphic to $\Fun[X_{1}, \ldots, X_{{\rm deg}(v)}]$.\medskip

Let us consider $v_i \in V(\Gamma)$, a vertex of degree $n_i$. Then, $\Spec(A_{v_i})$ is isomorphic to a $n_i$-dimensional affine space. We denote by $\Adj(v_i)$ the set of adjacent vertices of $v_i$, with cardinality $s_i$, by $E(v_i)$ the set of edges incident with $v_i$ and by $LE(v_i)$ the set of loose edges incident with $v_i$. Note that $s_i\leq n_i$ and that $LE(v_i)\subseteq E(v_i)$.  As $f$ is a graph automorphism, $f(v_i)$ is also a vertex $v_j$ of $\Gamma$ with degree $n_i$. We also know, as above, that $f(\Spec(A_{v_i}))=\Spec(A_{v_j})$ is isomorphic to an $n_i$-dimesional affine space and using the same terminology we consider the sets $E(v_j)$, $LE(v_j)$ and $\Adj(v_j)$ (also with cardinality $s_i$).\medskip

Then, $f$ induces a bijection between $E(v_i)$ and $E(v_j)$ and a bijection between $LE(v_i)$ and $LE(v_j)$. Since each edge incident with a vertex corresponds to a direction of the associated affine space, we can choose a base of each affine spaces $\Spec(A_{v_i})$ and $\Spec(A_{v_j})$ and $f$ will induce a unique bijection between both bases. We will call $\overline{f_i}$ the map

\begin{equation*}
\overline{f_i}: \Spec(A_{v_i}) \rightarrow \Spec(A_{v_j})
\end{equation*} 

\noindent induced by $f$ on the affine scheme associated to the vertex $v_i$. This map induces a morphism $\psi_i$ between the two corresponding $\Fun$-rings due to the contravariant functor between both categories:
\begin{equation*}
\psi_i: A_{v_j} \rightarrow A_{v_i}.
\end{equation*}

We know that the category of monoids ($\Fun$-rings) is a full subcategory of the category of sesquiads (see \cite{Deitmarcongruence}), so the morphism $\psi_i$ is also a morphism of sesquiads from $A_{v_j}$ to $A_{v_i}$ that induces in its turn a morphism of affine congruences schemes (\cite[theorem 3.2.1]{Deitmarcongruence})

\begin{equation*}
\overline{\psi_i}: \Spec_{c}(A_{v_i}) \rightarrow \Spec_{c}(A_{v_j})
\end{equation*}

\noindent in a contravariant way. So, we obtain for any vertex $v_k \in V(\Gamma)$ a homeomorphism $\overline{\psi_k}$ between affine congruence schemes and we can define a map $\widetilde{\psi}$ on $\mF(\Gamma)$ using the fact that, by definition, $\mF(\Gamma)$ is the union of the affine schemes associated to all vertices. Since all the morphisms $\overline{\psi_k}$ are morphisms of affine schemes and the intersection between the affine schemes are defined by the graph $\Gamma$, the union of morphisms is also well defined.\medskip

Once we have shown the functoriality property for {\em automorphisms} of a loose graph, we have to generalize the construction of $\widetilde{f}$ for a {\em general} morphism between two loose graphs. Let $\Gamma_1$ and $\Gamma_2$ be two loose graphs and let $\Hom(\Gamma_1, \Gamma_2)$ be the set of graph homomorphims from $\Gamma_1$ to $\Gamma_2$. Let us take $f \in \Hom(\Gamma_1, \Gamma_2)$ and use the same notation as above.\medskip

Notice that we may assume $f$ to be surjective. Otherwise we restrict to the image $f(\Gamma_1)\subset \Gamma_2$ and, using the (COV) property, we obtain the required map (composition with an embedding):

\begin{center}
\begin{tikzpicture}[auto]
\matrix(rings) [matrix of math nodes, row sep=0.6cm, column sep=0.4cm, ampersand replacement=\&]
{\mF(\Gamma_1)\& \& \mF(f(\Gamma_1))\\
\& \& \\
\& \& \mF(\Gamma_2)\\};
\draw[->] (rings-1-1) to node {\footnotesize{$\mF(f)$}} (rings-1-3);
\draw[dotted, ->] (rings-1-1) to [left=3pt] (rings-3-3);
\draw[right hook->] (rings-1-3) to node [right=3pt]{\footnotesize{$i$}} (rings-3-3);
\end{tikzpicture}
\end{center}

For general loose graph morphisms, the degree of vertices does not have to be preserved. When the degree of $v_i$ and $f(v_i)=v_j$ are equal, we define the morphism $\overline{f_j}$ as before. For the remaining case, let us denote by $S_{v}$ the {\em loose star} of the vertex $v$, i.e, the loose subgraph of $\Gamma_1$ formed by the vertex $v$ and all its incident edges. We consider then $f$ restricted to $S_{v}$, which corresponds to a morphism between two affine spaces on the scheme level. 

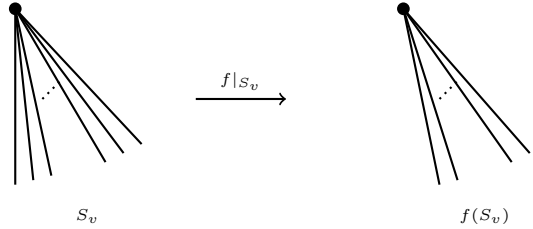
\begin{figure}[h]\label{p1}
\begin{tikzpicture}[style=thick, scale=1.2]
\draw (0,0.05)-- (0,2);
\draw (0.2,0.1)-- (0,2);
\draw (0.4,0.15)-- (0,2);
\draw (1,0.3)-- (0,2);
\draw (1.2,0.4)-- (0,2);
\draw (1.4,0.5)-- (0,2);
\draw [dotted] (0.3,1)-- (0.5, 1.2);
\fill (0,2) circle (2pt);
\draw [->] (2,1) -- (3,1);
\draw (4.7,0.05)-- (4.3,2);
\draw (4.9,0.1)-- (4.3,2);
\draw (5.5,0.3)-- (4.3,2);
\draw (5.7,0.4)-- (4.3,2);
\draw [dotted] (4.7,1) -- (4.9, 1.2);
\fill (4.3,2) circle (2pt);
\draw (2.5, 1.2) node {\tiny {$f|_{S_v}$}};
\draw (0.8, -0.3) node {\tiny {$S_v$}};
\draw (5.2, -0.3) node {\tiny {$f(S_v)$}};
\end{tikzpicture}
\end{figure}

Besides, $f(S_{v_i})$ is a proper loose subgraph of $S_{f(v_i)}$. Hence the following diagram gives us the desired homeomorphism:

\begin{center}
\begin{tikzpicture}[auto]
\matrix(rings) [matrix of math nodes, row sep=0.6cm, column sep=0.4cm, ampersand replacement=\&]
{\Spec(A_{v_i})\& \& \mF(f(S_{v_i}))\\
\& \& \\
\& \& \Spec(A_{v_j})\\};
\draw[->] (rings-1-1) to node {\footnotesize{$\mF(f|_{S_{v_i}})$}} (rings-1-3);
\draw[dotted, ->] (rings-1-1) to node [left=3pt]{\footnotesize{$\widetilde{f_i}$}} (rings-3-3);
\draw[right hook->] (rings-1-3) to node [right=3pt]{\footnotesize{$i$}} (rings-3-3);
\end{tikzpicture}
\end{center}

\medskip 
We reduced the study to the local restriction of $f$ to loose stars. Let us describe this situation in detail.
Suppose $v$ is a vertex of $\Gamma_1$ of degree $m$, and suppose the vertex of the loose star $f(S_v)$ has degree $n \leq m$ (this is always the case); then the morphism $f|_{S_{v_i}}$ is a loose graph morphism between the two loose stars $S_{v_i}$ and $f(S_{v_i})$.\medskip

We will call $\widetilde{f_v}$ the morphism induced between the corresponding affine spaces

\begin{equation*}
\widetilde{f_{v}} : \mathcal{F}(S_{v}) \longrightarrow \mathcal{F}(f(S_{v})).
\end{equation*} 

The morphism $\widetilde{f_v}$ is by definition a linear morphism between affine spaces with dimension the number of edges incident with the vertex of the corresponding loose star. Hence, we choose a basis $\{e_1, \ldots , e_m\}$ for $\mathcal{F}(S_{v}) \cong \mathbb{A}_{\Fun}^m$, where $e_i$ is the vector with 1 on the $i$-th coordinate and 0 elsewhere such that each element of the basis corresponds bijectively to an edge of the loose star $S_v$, and  $v$ corresponds to the point $(0, \ldots , 0)$. We do the same for the affine space $\mathcal{F}(f(S_{v})) \cong \mathbb{A}_{\Fun}^n$ and so we choose a basis $\{e'_1, \ldots, e'_n\}$.\medskip

Now that we have chosen a basis, we can easily describe the morphism $\widetilde{f_v}$ in terms of matrices. For each element $e_i$ of the basis of $\mathcal{F}(S_v)$, we consider the corresponding edge $g_i$ in $S_v$ and we set $\widetilde{f_v}(e_i)=e'_k$, where $e'_k$ is the element of the basis of $\mathcal{F}(f(S_v))$ associated to the edge $f(g_i)$. In the definition of loose graph morphism we allow contractions of edges having two end points, i.e., one edge with two end points might be contracted into the graph with one vertex. So it may happen that $f(g_i)$ is a vertex. But the only vertex existing on the loose star $f(S_v)$ is $f(v)$ so, in this case, we choose $\widetilde{f_v}(e_i)$ to be the zero vector.\medskip

Allowing contractions to be morphisms of loose graphs will let us have projections on the level of $\Fun$-schemes since, for instance, a projection of a projective line onto a point will be induced by the graph morphism sending the complete graph $K_2$ into one of vertex.\medskip

\begin{figure}[h]\label{p2}
\begin{tikzpicture}[style=thick, scale=1.2]
\draw (0,0)-- (2,0);
\fill (2,0) circle (2pt);
\fill (0,0) circle (2pt);
\draw [->] (3,0) -- (4,0);
\fill (5,0) circle (2pt);
\draw (3.5, 0.2) node {\tiny {$f|_{S_v}$}};
\end{tikzpicture}
\caption{Projection of $\mathbb{P}^1_{\Fun}$ on one point $P$.}
\end{figure}
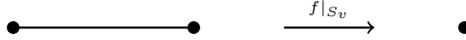

So locally $\widetilde{f_v}$ can be expressed by a matrix of size $n\times m$ whose columns are either the zero vector or a canonical vector, i.e., vector with only one non-zero entry. Reordering the basis $\{e_1, \ldots, e_m\}$, we obtain a block matrix of the form

\[ A_f:=\left(
\begin{array}{ccc|c|c|cccc|c}
0 & \cdots & 0 & A_1 & 0 & \cdots &  &\cdots & 0 & 0\\ \hline
0 &  \cdots & 0 & 0 & A_2 & 0 & &\cdots & 0 & 0\\ \hline
\vdots &&&\vdots & 0 &&&&& \vdots\\ 
 \vdots &&&& \vdots &&&&& \vdots\\ \hline
0 &\cdots & 0 & 0 & 0 & 0 &&\cdots & 0 & A_n\\ 
\end{array}
\right)
\]
\\

\noindent where the blocks $A_i$ are of size $1\times n_i$, with $n_i$ the number of vectors from $\{e_1, \ldots, e_n\}$ whose image is the vector $e'_i$, and have all 1-entries. Let us remark that $\sum^{n}_{j=1}n_j + n_0 = m$, where $n_0$ is the number of columns where all entries are 0. Note as well that if all edges are sent to edges by the morphism $f$, then $A_f$ has no zero part and if $n=m$, then $A_f$ is a nonsingular matrix. \medskip

These matrices $A_f$ are well defined over $\Fun$ since every column is a vector with at most one coordinate different from $0$. What is more, composition of two such morphisms corresponds to product of matrices. It is easy to verify that the product of two matrices of this form gives also a matrix with maximum one non-zero entry in each column, this entry being $1$, so composition of two morphisms is well defined.\medskip

After this construction, we need to remark that whenever the size of a block $A_i$ is bigger than $1\times 1$, then the image by such a linear morphism of a point of the affine space might be given by coordinates which include sums. Over $\Fun$ addition is not defined, but this problem is solved by considering the matrices defined over $\mathbb{F}_2$.  The reason why this consideration is also possible relies on the fact that points of a vector space over $\Fun$ (in the congruence setting of this paper) are exactly the same as points of the same vector space over $\F_2$. The difference between, e.g., projective spaces over $\Fun$ and over $\F_2$ can be seen geometrically on subvarieties of dim $\geq 1$ (and on the level of polynomial rings). \medskip

\subsection{Global action} The main idea to construct a morphism of schemes given a loose graph morphism $f$ from $\Gamma_1$ to $\Gamma_2$ is to define a morphism in their respective ambient spaces (the minimial projective spaces in which $\mathcal{F}(\Gamma_1)$ and $\mathcal{F}(\Gamma_2)$ are embedded) such that the restriction to the schemes $\mathcal{F}(\Gamma_1)$ and $\mathcal{F}(\Gamma_2)$ induces  a morphism between $\Fun$-schemes. Such a morphism will also induce the local mappings described in the previous subsection (just by considering their local action in the loose stars of vertices).\medskip

Consider now the completion $\overline{\Gamma_1}$, with $m_1+1$ vertices, and $\overline{\Gamma_2}$, with $m_2+1$ vertices, together with the embedding in their minimal projective space $\PG(m_1, \Fun)$ and $\PG(m_2, \Fun)$, respectively. Choose a set $\mathcal{R}_1 = \{e_0, \ldots, e_{m_1}, h \}$ of points of $\PG(m_1, \Fun)$ and a set $\mathcal{R}_2= \{e'_0, \ldots, e'_{m_2}, h' \}$ of points of $\PG(m_2, \Fun)$ in such a way that each $e_i$ and each $e'_j$ are the canonical vectors of $\PG(m_1, \Fun)$ and $\PG(m_2, \Fun)$, respectively, and each vertex of $\overline{\Gamma_1}$ corresponds to a canonical vector of $\mathcal{B}_1:=\mathcal{R}_1\setminus \{h\}$, and the same for the vertices of $\overline{\Gamma_2}$ and the vectors of $\mathcal{B}_2:=\mathcal{R}_2\setminus\{h'\}$. The elements $h$ and $h'$ are the ones having coordinates $[1 : \cdots : 1]$ w.r.t. their corresponding bases. Notice that when considering the extension of $\Fun$-schemes to $k$-schemes, the sets $\mathcal{R}_1$ and $\mathcal{R}_2$ would be skeletons of the projective spaces $\PG(m_1, k)$ and $\PG(m_2,k)$, respectively.\medskip

Let us describe the global construction. In this case we may also consider $f$ to be a surjective morphism without loss of generality, as before. The morphism $f$ induces a morphism $\overline{f}$ from $\overline{\Gamma_1}$ to $\overline{\Gamma_2}$. This morphism $\overline{f}$ sends every vertex of $\overline{\Gamma_1}$ to a vertex of $\overline{\Gamma_2}$. Besides, every element of $\mathcal{B}_1$ corresponds bijectively to a vertex of $\overline{\Gamma_1}$ and the same for elements of $\mathcal{B}_2$ and the vertices of $\overline{\Gamma_2}$, so reasoning as in the affine case and reordering the base $\mathcal{B}_1$, we get an $(m_2 + 1)\times (m_1 + 1)$-matrix of the form 

\[ P_f:=\left(
\begin{array}{c|c|cccc|c}
P_0 & 0 & \cdots &  &\cdots & 0 & 0\\ \hline
0 & P_1 & 0 & &\cdots & 0 & 0\\ \hline
\vdots & 0 &&&&& \vdots\\ 
& \vdots &&&&& \vdots\\ \hline
0 & 0 & 0 &&\cdots & 0 & P_{m_2}\\ 
\end{array}
\right)
\]
\\

Note that in this case, the matrix $P_f$ has no zero columns since all vertices are sent to vertices, i.e., every canonical vector of $\mathcal{B}_1$ is sent to a canonical vector of $\mathcal{B}_2$. It may also happen that en edge $e\in \overline{\Gamma_1}$ is contracted but this will imply that the two elements of $\mathcal{B}_1$ corresponding to the vertices of the edge are sent to the same element of the basis $\mathcal{B}_2$. As it happens for the affine case, the blocks $P_i$ are all of size $1\times n_i$, with all entries equal to 1 and $n_i$ being the number of elements of $\mathcal{B}_1$ whose image is the element $e_i'$. Besides, the identity $\sum^{n}_{j=0}n_j= m_1$ is also satisfied. (Nevertheless, in the general case when $f$ might not be surjective, $P_f$ might have zero rows, if an element of the basis $\mathcal{B}_2$ has no preimage.)\medskip

As in the affine case, those matrices are well defined over $\Fun$ and composition of morphisms, which is translated into a product of matrices, is also well defined. Notice that when you compose two morphisms it is not possible to reorder the frames so that both matrices are of the above form. Nevertheless, the product of two matrices having columns with only one non-zero element is also a matrix satisfying the same condition.\medskip

Hence, we obtain that for $f_1: \Gamma_1 \mapsto \Gamma_2$ and $f_2:\Gamma_2 \mapsto \Gamma_3$ two loose graph morphisms, the following property is satisfied:
\begin{equation*}
\mathcal{F}(g\circ f)=\mathcal{F}(g)\circ\mathcal{F}(f).
\end{equation*}

Let us remark that when we consider the schemes $\mathcal{X}_{1,k}=\mathcal{F}(\Gamma_1)_k$ and $\mathcal{X}_{2,k}=\mathcal{F}(\Gamma_2)_k$ over a field $k$ (or $\Z$), the morphism defined by the matrix $P_f$ also induces an action on the level of schemes over $k$ (or $\Z$). However, if we {\em first} consider the extension of $\Fun$-schemes to schemes over a field $k$ (or $\Z$) and then define a matrix $P_{f,k}$ in the same way as we did with $P_f$, one realizes that there exist many choices for $P_{f,k}$ inducing the same action on the basis vectors. So, in general, the following diagram

\begin{center}
\begin{tikzpicture}[auto]
\matrix(rings) [matrix of math nodes, row sep=0.6cm, column sep=0.4cm, ampersand replacement=\&]
{\mathcal{F}(\Gamma_1) \& \& \mathcal{F}(\Gamma_2)\\
\& \& \\
\mathcal{F}(\Gamma_1)_k \& \& \mathcal{F}(\Gamma_2)_k\\};
\draw[->] (rings-1-1) to node {\footnotesize{$P_f$}} (rings-1-3);
\draw[->] (rings-1-1) to node {\footnotesize{$\otimes k$}} (rings-3-1);
\draw[->] (rings-1-3) to node {\footnotesize{$\otimes k$}} (rings-3-3);
\draw[->] (rings-3-1) to node {\footnotesize{$ P_{f,k} $}} (rings-3-3);

\end{tikzpicture}
\end{center}

\noindent is not commutative. One way of dealing with this problem is to define a morphism on the level of $k/\mathbb{Z}$-schemes as the class $[P_{f,k}]$ of morphisms having the same action on the basis vectors. As such, we obtain a well-defined functor from the category of loose graphs to the category of $k/\mathbb{Z}$-schemes making the previous diagram commutative.

\subsection{Different categories for projective spaces} 

After the previous construction one could realize that these morphisms might not be injective on the level of projective spaces so we need to choose the concrete category of projective spaces that we want to work with. We will see in this subsection how one should interpret the morphisms depending on the chosen category. We only work over $\Fun$; similar considerations over ``real fields'' follow easily (and are in fact {\em easier}). 

\subsubsection*{Category with injective linear maps} We consider the category of projective spaces whose morphisms are injective linear maps. In this case $P_f$ defines an injective linear map if and only if $m_2 \geq m_1$ and the rank of $P_f$ equals $m_1 + 1$ over $\mathbb{F}_2$ (since it is a homogeneous linear system of equations). These conditions is equivalent to say that every element of the basis $\mathcal{B}_1$ is sent to a different element of the basis $\mathcal{B}_2$ and, by the bijection described above between basis and vertices of the graphs, we have that $P_f$ is induced from an injective morphism of graphs. \medskip

So in this case, our functor $\mathcal{F}$ will be a functor between the category of graphs with injective morphisms and the category of congruence schemes with injective linear morphisms.

\subsubsection*{Category with rational maps} 

In the second case we consider the category of projective spaces whose morphisms are {\em rational maps}. A {\em rational map} $f: V \rightarrow W$ between two varieties is an equivalence class of pairs $(f_U, U)$ in which $f_U$ is a morphism of varieties defined from an open set of $U\subseteq V$ to $W$, and two pairs $(f_U, U)$ and $(f_{U'}, U')$ are equivalent if $f_{U}$ and $f_{U'}$ coincide on $U\cap U'$. We adapt the same nomenclature for other types of schemes, such as Deitmar schemes (with possible congruences).\medskip

Consider now the linear map $P_f$ defined in the previous subsection and suppose it has a nontrivial kernel. One should remember that the map $P_f$ is defined in the minimal projective spaces in which $\mathcal{F}(\Gamma_1)$ and $\mathcal{F}(\Gamma_2)$ are embedded (denoted $\PG(m_1, \Fun)$ and $\PG(m_2, \Fun)$, respectively). That implies we only have to consider the case where the kernel of $P_f$ intersects with the scheme $\mathcal{F}(\Gamma_1)$, since in any other cases the induced map on the schemes will be injective.\medskip

We will then prove that the kernel of $P_f$ is a closed subset of  $\PG(m_1, \Fun)$. By relative topology its intersection with $\mathcal{F}(\Gamma_1)$ will be closed in $\mathcal{F}(\Gamma_1)$. We define the kernel of $P_f$ over $\Fun$, which we denote by $\mathrm{ker}(P_{f})_{\Fun}$, as

\begin{equation*}
\mathrm{ker}(P_{f})_{\Fun}:=\{x\in \mathcal{F}(\Gamma_1)\otimes \mathbb{F}_2 ~|~ x\in \mathrm{ker}(P_f) \}.
\end{equation*}

Notice that to define the map $P_f$ one has to consider the matrix defined over $\mathbb{F}_2$ and that, in terms of points, there is a bijective correspondence between the points of the scheme $\mathcal{F}(\Gamma_1)\otimes_{\Fun}\mathbb{F}_2$ and the points of the congruence scheme $\mathcal{F}(\Gamma_1)$. So in this way the kernel of $P_f$ over $\Fun$ is well defined. \medskip

To prove that $\mathrm{ker}(P_{f})_{\Fun}$ is a closed subset of $\mathcal{F}(\Gamma_1)$ we will verify that every point of the set $\mathrm{ker}(P_{f})_{\Fun}$ is indeed closed in the congruence scheme $\mathcal{F}(\Gamma_1)$. For, take a point $x\in \mathrm{ker}(P_{f})_{\Fun}$. Considering $x$ as a point in the projective space gives us its coordinates; let us write $x=[a_0: \cdots :a_{m_1}]$ (not all entries 0) and let $a_{i_0}$ be the first coordinate equal to $1$.
Then $x$ defines a congruence $\mathcal{C}_x$ in the projective congruence scheme corresponding to $\PG(m_1, \Fun)$, given by $\langle x_{i}\sim 0 \mbox{ if } a_i=0, ~ x_{i}\sim x_{i_0}\mbox{ if } a_i=1 \rangle$. It is a homogeneous maximal congruence in the Zariski topology (remember that a maximal congruence in projective schemes is maximal w.r.t. ``not containing the irreducible congruence''). (To see the construction of congruence projective schemes, we refer to section $2.3$.)\medskip

We have proved that a point $x\in \mathrm{ker}(P_{f})_{\Fun}$  is a closed point in the projective scheme  $\PG(m_1, \Fun)$, so the set $\mathrm{ker}(P_{f})_{\Fun}$ is closed as well since it is a finite union of closed sets.
So if $f$ is a morphism of loose graphs,  $\mathcal{F}(f)=P_f$ is a rational map defined on $U=\mathcal{F}(\Gamma_1)\setminus \mathrm{ker}(P_{f})_{\Fun} $. So, the functor $\mathcal{F}$ will be a functor between the category of loose graphs with morphisms of loose graphs and the category of congruence schemes with rational maps.

\section{Different types of automorphisms}

\subsection{Projective automorphism group}
Let $\Gamma$ be a loose graph and $\mathcal{F}(\Gamma)$ be its $\Fun$-scheme.
We define the {\em projective automorphism group} of the scheme $\mathcal{X}_k$, denoted by $\mbox{Aut}^{\mbox{\tiny proj}}(\mathcal{X}_k)$, as the group of automorphisms of the ambient projective space of $\mathcal{X}_k$ stabilizing $\mathcal{X}_k$ setwise, modulo the group of such automorphisms acting trivially on $\mathcal{X}_k$.

\subsection{Combinatorial automorphism group}
Let $\Gamma$ be a loose graph and consider $\mathcal{F}(\Gamma)$. 
To define the {\em combinatorial automorphism group} of the scheme $\mathcal{X}_k=\mathcal{F}(\Gamma)\otimes_{\Fun} k$, with $k$ a field, we want to consider the scheme as an {\em incidence geometry of rank 2}, i.e., as a set of points $\mathcal{P}$ and a set of lines $\mathcal{L}$ in which a relation of incidence is given.
We consider the set of points to be the set of $k$-rational points of $\mathcal{X}_k$ and the set of lines to be consisting of both projective lines (over $k$) and {\em complete affine lines}. A {\em complete affine} line $l$ of $\mathcal{X}_k$ is a line whose projective completion $\bar{l}$ intersects the scheme $\mathcal{X}_k$ in the whole projective line $\bar{l}$ minus one point.\medskip

A {\em combinatorial automorphism} $g$ of $\mathcal{X}_k$ is a bijective map on the set of points and on the set of lines preserving incidence, i.e., if a point $p$ is on the line $L$ of $\mathcal{X}_k$, then the image of $p$ is on the image of $L$ in $\mathcal{X}_k$. We will denote by $\mbox{Aut}^{\mbox{\tiny comb}}(\mathcal{X}_k)$ the group of combinatorial automorphisms of the scheme $\mathcal{X}_k$.\medskip

The next two results show that combinatorial automorphisms automatically preserve the linear subspace structure of the schemes.

\begin{observ}\label{obs1} {\em Let $\mathcal{X}_k$ be a scheme coming from a loose tree and let $g$ be a combinatorial automorphism of $\mX_k$. If $\A$ is a $d$-dimensional affine space contained in $\mathcal{X}_k$, then $\A^g$ is also an affine space of dimension $d$, contained in $\mathcal{X}_k$ and isomorphic to $\A$.} 
\end{observ}

\prf To prove that $\A^g$ is an affine space isomorphic to $\A$, it is sufficient to recall the axiomatic definition of an affine space in terms of an incidence geometry of rank 2 in which one has a set of points $\mathcal{P}$, a set of lines $\mathcal{L}$ and an equivalence relation ``$\|$'' of parallelism defined on 
the set of lines. The idea is that using the axioms one observes that $\A^g$ is an axiomatic affine space with the same dimension of $\A$. It is then obvious that $g$ is an isomorphism between (axiomatic) affine spaces (and so $\A^{g}$ is also defined over $k$). The axioms are the following:

\begin{itemize}
	\item Each pair $P,Q$ of distinct points is contained in a unique line $l$.
	\item  For each point $P$ and each line $l$, there is a unique line $l'$ such that $P\in l'$ and $l \| l'$.
	\item {\em Trapezoid axiom}. Let $\overline{PQ}$ and $\overline{RS}$ be distinct parallel lines and let $T$ be a point of $\overline{PR}\setminus\{P,R\}$. Then, there must be a point incident with $\overline{PQ}$ and $\overline{TS}$.
	\item {\em Parallelogram axiom}. If no line has more than two points, and if $P$, $Q$ and $R$ are three distinct points, then the line through $R$ parallel to $PQ$ must have a point in common with the line through $P$ parallel to $QR$.
	\item {\em Thickness.} Each line contains at least two points.
	\item {\em Space axiom}. There exists two disjoint lines $l$ and $l'$ such that $l\nparallel l'$. Notice that this axiom is only required if the dimension of the affine space is greater than 2.
\end{itemize}

Every axiom is satisfied in $\A^g$ since the automorphism $g$ preserves the incidence relations and $g$ is injective on points and lines of $\A$. Hence, $\A^g$ is an affine space. It remains to prove that it is indeed of dimension $d$.\medskip

Recall that the {\em geometric dimension} of the affine space $\A$ is given recursively by the largest number ($d$ in this case) for which there exists a strictly ascending chain of subspaces of the form:

\begin{equation}\label{paraA}
\emptyset \subset X_0 \subset X_1 \subset \cdots \subset X_d=\A,
\end{equation}

\noindent where $X_i$ is a subspace of geometric dimension $i$. Since $g$ is an automorphism, applying $g$ to this chain we will obtain a new chain of the form 

\begin{equation*}
\emptyset \subset X_0^g \subset X_1^g \subset \cdots \subset X_d^g=\A^g,
\end{equation*}

\noindent where all subspaces $X^g_i$ are of dimension greater or equal to $i$, since $g$ is injective. Let us now suppose that dim$(A^g)=j>d$; then there will exist a chain of the form 

\begin{equation*}
\emptyset \subset Y_0 \subset Y_1 \subset \cdots \subset Y_j=\A^g.
\end{equation*}

Applying $g^{-1}$ to this new chain, we will obtain a chain for $\A$ longer than (\ref{paraA}) since $g^{-1}$ is injective as well. But this is not possible since (\ref{paraA}) is a chain of maximal length. 

\begin{observ}\label{obs2}{\em With the same conditions of the previous observation, if $\hP$ is a $d$-dimensional projective space contained in $\mathcal{X}_k$, then $\hP^g$ is also a projective space of dimension $d$, contained in $\mathcal{X}_k$ and isomorphic to $\hP$.}
\end{observ}

\prf As before we just have to recall the axiomatic definition of a projective space in terms of an incidence geometry of rank 2. The axioms in this case are the following:

\begin{itemize}
	\item Two different points are exactly in one line.
	\item {\em Veblen's axiom}. If $a$, $b$, $c$ and $d$ are different points and the lines $ab$ and $cd$ meet, then so do the lines $ac$ and $bd$.
	\item {\em Thickness}. A line has at least 3 points.
\end{itemize}

For the same reason as in the affine case, every axiom is satisfied in $\hP^g$, so $\hP^g$ is a projective space. The fact that dim$(\hP^g)=d$ is proven in the same way as for the affine case.\eop\medskip

Let $r:=\mbox{max}\{ \mbox{deg}(v) ~|~ v\in V(\Gamma)$ and $v$ defines an affine space $\A_v$ sucht that $\overline{\A_v}$ is not contained in $\mathcal{X}_k\}$ and $s:=\mbox{max}\{n-1 ~|~ K_n\subseteq\Gamma\}$. We will consider $\mathcal{X}_k$ as an incidence geometry of ``{\em double rank} $(r,s)$.'' We define $\mathcal{X}_k$ as the $(r+s+2)$-tuple $(K, A_1, \ldots , A_{r}, P_1, \ldots , P_{s}, {\rm \bf I} )$, where $A_i$ is the set of $i$-dimensional affine subspaces of $\mathcal{X}_k$ whose completion is not contained in $\mathcal{X}_k$, $P_k$ is the set of $k$-dimensional projective subspaces of $\mathcal{X}_k$, $K=A_0=P_0$, and ${\rm \bf I}$ is the natural incidence relation between these spaces. Note that the sets $A_i$ and $P_j$ are non empty for all $i,j$. \medskip

If $\mathcal{X}_k$ is, e.g., a projective space of dimension $d$, then the double rank is $(0,d)$. If $\Gamma$ is a tree, then the double rank is $(r,1)$ or $(0,0)$ (if $\Gamma$ is a vertex).\medskip

With this definition of $\mathcal{X}_k$ the two previous observations lead to the following result.

\begin{corollary} Let $\Gamma$ be a loose graph, $\mathcal{X}_k$ its corresponding scheme over $k$ defined by $\mathcal{F}(\Gamma)\otimes_{\Fun} k$ and $g$ a combinatorial automorphism of $\mathcal{X}_k$. If we define the numbers $r$ and $s$ as above, then $g$ is also an automorphism of $\mathcal{X}_k$ as an incidence geometry of double rank $(r,s)$. \end{corollary}

\prf The proof of this corollary follows immediately after Obervation \ref{obs1} and Observation \ref{obs2} and the fact that $g$ preserves incidence relations when $\mX_k$ is considered as an incidence geometry of rank $2$.

\subsection{Topological automorphism group}

We define a {\em topological automorphism} $g$ of the scheme $\mathcal{X}_k$ as a homeomorphism of its underlying topological space, i.e, a bijective continuous map with a continuous inverse map. In a natural way, we obtain the {\em topological automorphism group} of $\mathcal{X}_k$, denoted by $\mbox{Aut}^{\mbox{\tiny top}}(\mathcal{X}_k)$.\medskip

First note that for any field $k$, the closed set topology of $\Spec(k[X])$ consists of $(0)$ and all closed points (since all prime ideals are maximal) and all finite sets of such points that contain $(0)$. So $\Aut^{\text{top}}(\Spec(k[X]))$ is isomorphic to the symmetric group on the set of closed points. On the other hand, $\Aut^{\text{comb}}(\Spec(k[X])$ is isomorphic to the symmetric group on the $k$-rational points, so as soon as $k$ is not algebraically closed, the groups are not the same. Now let $\mX_k$ be $\Spec(k[X_1,\ldots,X_m])$ with $m \geq 2$, and let $U$ be an affine subline. Then $\Aut^{\text{top}}(\mX_k)$ induces $\Aut^{\text{top}}(U)$ on the topology of $U$, which, as we have seen, is isomorphic to the symmetric group on the closed points of $U$. The combinatorial automorphism group of $\mX_k$ induces
the affine group $\mathbf{A\Gamma L}_1(k)$ on $U$ (acting on the $k$-rational points). So in general the groups are not isomorphic.

The next proposition deals with the other direction.

\begin{proposition} 
\label{combtopo}
The combinatorial group of a scheme $\mX_k$ is a subgroup of the topological automorphism group of $\mathcal{X}_k$. 
\end{proposition}\medskip

\prf 
Let us first take a combinatorial automorphism $f$ of $\mX_k$. We can reduce our proof w.l.o.g. to the case of an affine space defined by one of the loose stars corresponding to a vertex of $\Gamma$. This is possible since an automorphism of the scheme $\mX_k$ can be constructed as the union of the local morphisms of affine spaces. \medskip

Let $\Spec(A_{v})$ be the affine space corresponding to the vertex $v\in \Gamma$. If $f$ is a combinatorial automorphism of $\mX_k$, by (\ref{obs1}) we know that $f$ induces also a combinatorial isomorphism $f_v$ from $\Spec(A_v)\otimes_{\Fun} k$ to $\Spec(A_{f(v)})\otimes_{\Fun} k$. Remember that due to the definition of the functor $\mathcal{F}$ the affine spaces $\Spec(A_v) \otimes_{\Fun} k$ and $\Spec(A_{f(v)})\otimes_{\Fun} k$ are isomorphic to $\Spec(k[X_1, \ldots , X_n])$ and $\Spec(k[Y_1, \ldots , Y_n])$, with $n=\mbox{deg}(v)=\mbox{deg}(f(v))$, respectively. \medskip
Hence, the latter combinatorial isomorphism $f_v$ induces a ring isomorphism between the corresponding coordinate rings $k[X_1, \ldots ,X_n]$ and $k[Y_1,\ldots, Y_n]$ (induced by the action on the coordinate hyperplanes) that gives, by functoriality, an isomorphism of affine schemes from $\Spec(A_v)\otimes_{\Fun} k$ to $\Spec(A_{f(v)})\otimes_{\Fun} k$. For each vertex $v$ of the graph $\Gamma$ we hence obtain an induced topological isomorphism between the local affine $k$-schemes corresponding to $v$ and $f(v)$. By considering the union of these isomorphisms we finally obtain the topological automorphism of the scheme $\mX_k$.\eop \\


\medskip
\section{Toy example}

Let $\Gamma$ be the connected loose graph on two vertices ($x$ and $y$) of regular degree $2$. In this section we show that 
\begin{equation}
\Aut(\mF(\Gamma) \times_{\Fun}k) \cong \Aut^{\mathrm{proj}}(\mF(\Gamma) \times_{\Fun}k)
\end{equation}
for any field $k$. Here (and throughout), $\Aut(\cdot)$ is the {\em combinatorial} automorphism group.

For the rest of this section, fix a field $k$.
Denote the loose edge on $x$ by $L_x$, the loose edge on $y$ by $L_y$, and $xy$ by $L$. Also, let $\A_x$ and $\A_y$ be the affine planes
corresponding (respectively) to  the vertices $x$ and $y$ through $\mF$. To shorten notation, we will write $\mX$ for $\mF(\Gamma)$, and $\mX_k$
for $\mF(\Gamma) \times_{\Fun} k$ (so that in particular $\mX_{\Fun} = \mX$). 

For now, we want to see $\mX_k$ coming together with its embedding 
\begin{equation}
\mX_k \ \ \hookrightarrow\ \ \PG(3,k).
\end{equation}

It makes sense to projectively complete $\A_x$ and $\A_y$ | that is, to add the respective lines at infinity $X$ and $Y$; obviously, any 
element of $\Aut^{\mathrm{proj}}(\mX_k)$ also fixes the ``projective completion'' $\overline{\mX_k}$. As $X$ is incident with $y$ (as a point of $\mX_k$) and $Y$ with $x$ (as a point of $\mX_k$), it is now clear that $\alpha \in \Aut^{\mathrm{proj}}(\mX_k)$ (with $\alpha \in \mathbf{P\Gamma L}_4(k)$) if and only if 
$\alpha$ stabilizes the configuration $(\{x,y\},\{X,Y\},\{(x,Y),(Y,x),(y,X),(X,y) \})$ $:= \rho$ (defined as incidence geometry). Call such a configuration a {\em root}, and denote it also by $(Y,x,xy,y,X)$. It is very important to notice that $X$ and $Y$ are {\em projective lines}, and not {\em affine} lines.

In this paper, if $\fP$ is a projective space, and $\Omega$ is a hyperplane, by $T(\Omega)$ we denote the group of translations of $\fP$ with axis $\Omega$; it fixes $\Omega$ pointwise, acts sharply transitively on the points of $\fP \setminus \Omega$, and is a subgroup of $\PGL(\fP)$.

\begin{proposition}
\label{transroot}
$\mathbf{P\Gamma L}_4(k)$ acts transitively on the roots of $\PG(3,k)$.
\end{proposition}
{\em Proof}.\quad
Obviously $\mathbf{P\Gamma L}_4(k)$ acts transitively on the ordered triples $(u,uv,v)$, with $u \ne v$ points of $\PG(3,k)$ (as it acts transitively on the lines, and a line stabilizer induces the natural action of $\mathbf{P\Gamma L}_2(k)$, which is $3$-transitive).
Fix such a triple $(x,xy,y)$. Let $Y, Y'$ be different lines on $x$, both different from $xy$.  Consider a plane $\nu$ containing $xy$ but not $Y$ nor $Y'$. Then there is an element in $T(\nu)$ that maps $Y'$ to $Y$, so from now on, we also fix $Y$.  Now let $X, X'$ be lines on $y$ different from $xy$, and not meeting $Y$. Define the plane $\rho := \langle Y,xy \rangle$, and note that it does not contain $X$ nor $X'$. Then $T(\rho)$ contains an element which maps $X'$ to $X$. The claim follows.
\eop \\

Note that roots are {\em ordered}.

By the proof of the previous proposition, we immediately have the following.

\begin{corollary}
$\mathbf{PGL}_4(k)$ acts transitively on the roots of $\PG(3,k)$.
\end{corollary}
{\em Proof}.\quad
One can replace $\mathbf{P\Gamma L}_4(k)$ by $\PGL_4(k)$ in the proof of Proposition \ref{transroot}. Furthermore, all translations are elements in $\PGL_4(k)$.
\eop \\

The following is immediate.

\begin{proposition}
The kernel of the action of $\mathbf{P\Gamma L}_4(k)_{\mX_k}$ on $\mX_k$ is trivial.\eop
\end{proposition}

\begin{proposition}
Let $\hP_x$ be the projective $k$-plane generated by $x, xy$ and $Y$. (For later purposes, we similarly define $\hP_y$.) 
Let $A := \Aut(\hP_x)_{(Y,x,xy,y)}$ be the elementwise stabilizer of $\{x,y,xy,Y\}$ in $\Aut(\hP_x)$, where the latter is the incidence geometrical (= combinatorial) automorphism group of $\hP_x$ (so isomorphic to $\mathbf{P\Gamma L}_3(k)$). 
Then each element of $A$ extends to an element of $\PGaL_4(k)_{\mX_k}$ (in a not necessarily unique fashion).
\end{proposition}
{\em Proof}.\quad
Let $\alpha \in A$ be arbitrary; then $\alpha$ extends to  elements of $\PGaL_4(k)$, for instance to $\widetilde{\alpha}$. Note that $\widetilde{\alpha}$ fixes $y$. Suppose that $X^{\widetilde{\alpha}} =: X'$. Now let $\beta$ be an element in $T(\hP_x)$ which maps $X'$ back to $X$; then $\beta \circ \widetilde{\alpha}$ fixes the root $(Y,x,xy,y,X)$ and induces $\alpha$ on $\hP_x$.
\eop \\

\begin{remark}{\rm
It is important to note that $\Aut(\hP_x)$ coincides  with the automorphism group of $\hP_x$ induced by the automorphisms of $\mathbf{P \Gamma L}_4(k)$.
}\end{remark}

The ``number'' of ways to extend an element $\alpha$  is easy to determine. For, if $\gamma$ and $\gamma'$ are two such elements, then $\gamma^{-1} \circ \gamma'$ fixes $\hP_x$ pointwise, while fixing $X$. This group faithfully induces $\PGL_2(k)_y$ on the projective line $X$ (we are in the projective group, since $\hP_x$ is pointwise fixed). Its order is $\vert k \vert(\vert k \vert - 1)$ (the group acts sharply $2$-transitively on $X \setminus \{y\}$ and is isomorphic to $k \rtimes k^{\times}$). 

We now have all the ingredients for writing down $\Aut^{\mathrm{proj}}(\mX_k)$. First of all, it is clear that  $\Aut(\Gamma) \cong \langle \varphi \rangle$, with $\varphi \ne \id$ an  involution.  By Proposition \ref{transroot}, there is an element in $\PGaL_4(k)$ which stabilizes the root $\rho$, and which has the same action as $\varphi$. And obviously, the subgroup of $\Aut^{\mathrm{proj}}(\mX_k)$ which fixes both $x$ and $y$ is a normal subgroup of $\Aut^{\mathrm{proj}}(\mX_k)$.

\begin{theorem}
\label{DDC}
Let $\PGaL_2(k)$ be the automorphism group of the projective line $\PG(1,k)$, and let $u, v$ be distinct points of the latter.
Let $C := \PGaL_2(k)_{(u,v)} \cong k^{\times} \rtimes \Aut(k)$, and $D := \PGL_2(k)_u \cong k \rtimes k^{\times}$. Then 
\begin{equation}
 \Aut^{\mathrm{proj}}(\mF(\Gamma) \times_{\Fun}k) \cong (D \rtimes (D \rtimes C)) \rtimes \langle \varphi \rangle. 
\end{equation}
In the latter expression,
\begin{itemize}
\item
$(D \rtimes (D \rtimes C))$ is the elementwise stabilizer of $\{x,y\}$ in $ \Aut^{\mathrm{proj}}(\mX_k)$;
\item
the ``first $D$'' is $T(\hP_x) \cap \Aut^{\mathrm{proj}}(\mX_k)$;
\item
$E := D \rtimes C$ is the pointwise stabilizer of $X$ in  $\Aut^{\mathrm{proj}}(\mX_k)$;
\item
$D$ (``in $E$'') is the poinwtise stabilizer of $xy$ in $E$, and $C$ (``in $E$'') is the action induced by $E$ on $xy$.
\end{itemize}
\eop
\end{theorem}

For later purposes, we need an approach which allows a possibility to extend to more general cases.
Let $\alpha \in \Aut^{\mathrm{proj}}(\mX_k)$; then $\alpha$ induces an element $\alpha_x$ of $\Aut(\hP_x)$ which fixes $(Y,x,xy,y)$, and also an 
element $\alpha_y$ of $\Aut(\hP_y)$ which fixes $(x,xy,y,X)$, and both elements have the same action on the projective line $xy$. And vice versa, we have that  $\Aut^{\mathrm{proj}}(\mX_k)$ is completely determined by the data

\begin{equation}
\label{equi}
\Big\{ (\alpha_x,\alpha_y)\ \vert \ \alpha_x \in \Aut(\hP_x)_{(Y,x,xy,y)}, \alpha_y \in \Aut(\hP_y)_{(x,xy,y,X)}, {\alpha_x}_{\Big\vert xy} \equiv {\alpha_y}_{\Big\vert xy} \Big\}.
\end{equation}

Before using this observation, we prove the next theorem.

\begin{theorem}\label{6.7}
Let $\Gamma$ be the connected loose graph on two vertices ($x$ and $y$) of regular degree $2$. Then
\begin{equation}
\Aut(\mF(\Gamma) \times_{\Fun}k) \cong \Aut^{\mathrm{proj}}(\mF(\Gamma) \times_{\Fun}k)
\end{equation}
for any field $k$.
\end{theorem}
{\em Proof}.\quad
Recall that it is obvious by mere definition that $\Aut^{\mathrm{proj}}(\mX_k) \leq \Aut(\mX_k)$.
Let $\gamma$ be an element of $\Aut(\mX_k) \setminus \Aut^{\mathrm{proj}}(\mX_k)$; then there also exists an element $\gamma'$ in $\Aut(\mX_k) \setminus \Aut^{\mathrm{proj}}(\mX_k)$ which fixes both $x$ and $y$ (that is, which fixes the root $(Y,x,xy,y,X)$ elementwise). For, it is obvious that there is an $\epsilon \in \Aut^{\mathrm{proj}}(\mX_k)$ which switches $x$ and $y$ (and $Y$ and $X$) (by using Proposition \ref{transroot}). If $\gamma$ already fixes $x$ and $y$, there is nothing to prove. If not, $\epsilon \circ \gamma$ fixes $x, y$, and is not in $\Aut^{\mathrm{proj}}(\mX_k)$.
Now $\gamma'$ induces an element $\gamma'_x$ in  $\Aut(\hP_x)_{(Y,x,xy,y)}$ and an element $\gamma'_y$ in  $\Aut(\hP_y)_{(x,xy,y,X)}$ which agree on $xy$. We have seen that there exists an element $\gamma^*$ in 
$\Aut^{\mathrm{proj}}(\mX_k)$ which also yields the data $(\gamma'_x,\gamma'_y)$; composing $\gamma'$ with ${\gamma^*}^{-1}$, we obtain the identity of $\Aut(\mX_k)$. The isomorphism follows.
\eop \\

If $L$ is a line of $\PG(3,k)$, by $\PGaL_4(k)_{[L]}$ we will denote the pointwise stabilizer of $L$ in $\PGaL_4(k)$. (Note that it is a subgroup of 
$\PGL_4(k)$.) More generally, if $S$ is a set of points in $\PG(3,k)$, $\PGaL_4(k)_{[S]}$ denotes its pointwise stabilizer (and this is not necessarily a subgroup of $\PGL_4(k)$).

\begin{lemma}
Define $A := \Aut^{\mathrm{proj}}(\mX_k) \cap \PGaL_4(k)_{[Y]}$, and  $B := \Aut^{\mathrm{proj}}(\mX_k) \cap \PGaL_4(k)_{[X]}$.
Then $A \cong \Aut(\hP_y) \cap \PGL_3(k)$ (where it is obvious what we mean by the latter expression, namely the projective general elements in $\Aut(\hP_y)$), and   $B \cong \Aut(\hP_x) \cap \PGL_3(k)$.
\end{lemma}
{\em Proof}.\quad
We prove the assertion for $A$. Let $\alpha$ be any element in $\Aut(\hP_y) \cap \PGL_3(k)$; we have seen that $\alpha$ extends to some element $\widetilde{\alpha}$ of $\Aut(\mX_k)$, and that any such element induces a projective general linear element on $Y$. So  there is a unique element in
$\PGaL_4(k)_{[\hP_y]}$ with the same action on $X$. Composing with the inverse of $\widetilde{\alpha}$, we obtain an element of $A$ which induces $\alpha$ on $\hP_y$. The required isomorphism easily follows.  
\eop \\

\begin{theorem}
\label{thmcp}
Let $\PGL(\mX_k)_{(x,y)}$ be defined as 
\begin{equation}
\Aut^{\mathrm{proj}}(\mX_k)_{(x,y)} \cap \PGL_4(k).
\end{equation}
Then $\PGL(\mX_k)_{(x,y)}$ is isomorphic to the internal central product of $A$ and $B$. 
\end{theorem}
{\em Proof}.\quad
It is obvious that $\langle A, B \rangle = \PGL(\mX_k)_{(x,y)}$, so we only have to show that $[A,B] = \{\id\}$. Now if $a \in A$ and $b \in B$, we have that 
$[a,b] = a^{-1}b^{-1}ab$ fixes $Y$ and $X$ pointwise. On the other hand, both $a$ and $b$ induce elements in $\Aut(xy)_{(x,y)} \cap \PGL_2(k) \cong k^{\times}$, and this is an abelian group. So $[a,b]$ acts as the identity on $xy$. It now easily follows that $[a,b]$ acts trivially on $\mX_k$.
\eop \\

In general, we have the next conclusion.

\begin{theorem}
We have that
\begin{equation}
\Aut^{\mathrm{proj}}(\mX_k)\ \cong\ \Aut(\mX_k)\ \cong\ ((A * B) \rtimes \Aut(k)) \rtimes \langle \varphi \rangle. 
\end{equation}
\end{theorem}
{\em Proof}.\quad
Follows from Theorem \ref{DDC}, and the identities
\begin{equation}
\PGL(\mX_k)_{(x,y)}\ \unlhd\  \Aut^{\mathrm{proj}}(\mX_k)_{(x,y)}\ \unlhd\ \Aut^{\mathrm{proj}}(\mX_k).
\end{equation}
\eop

\medskip
\section{Trees}

If $\Gamma$ is a connected loose tree, and $k$ a field, one of the first things to hope is that:
\begin{itemize}
\item
$\Aut(\mX_k)$ {\em acts} on the set of affine spaces defined by the vertices $\Gamma$;
\item
this action is induced by $\Aut(\Gamma)$.
\end{itemize}

These properties are not true in general | look for instance at a projective plane (coming from a triangle): for no field $k \ne \Fun$ one has that $\Aut(\mX_k)$
induces an action on the three subplanes corresponding to the vertices.

If the toy example generalizes naturally, one candidate for $\Aut^{\mathrm{proj}}(\mX_k)$ would be
\begin{equation}
(U \rtimes \Aut(k)) \rtimes \Aut(\Gamma)^*,
\end{equation}
where $U \rtimes \Aut(k)$ is the part that fixes all vertices of $\Gamma$ (once pulled to $k$), and $U$ is the projective general linear part of the latter.
After the toy example, $U$ should be isomorphic to a central product of the appropriate groups. Also, $\Aut(\Gamma)^*$ is the automorphism group of the graph underlying $\Gamma$.\\

The first thing to do is generalize the little theory of roots.

\subsection{Fundaments}

Consider a $\PG(a + b - 1,k) = \pi$ over the field $k$, with $a, b \geq 2$. A {\em fundament} of {\em type} $(a,b)$ of $\pi$ is a triple $(\alpha,xy,\beta)$, where $\alpha$ is an $(a - 1)$-dimensional projective subspace of $\pi$, $\beta$ a $(b - 1)$-dimensional projective subspace, and $xy$ a projective line for which $\alpha \cap xy = \{x\}$ and 
$\beta \cap xy = \{y\}$, and
such that 
\begin{equation}
\langle \alpha,xy \rangle \cap \langle \beta,xy \rangle = xy.
\end{equation}

Note that $\langle \alpha,\beta \rangle = \pi$, and that a fundament of $\PG(3,k)$ is a root. A fundament {\em with ends} is a $5$-tuple $(\alpha,A,xy,\beta,B)$
where $(\alpha,xy,\beta)$ is a fundament (of type $(a,b)$), $A$ a projective subspace  of $\alpha$ which does not contain $x$, and $B$ is a projective subspace of $\beta$ not containing $y$. Such a fundament has {\em type} $(a,b;c,d)$ if $A$ and $B$ respectively have dimension $c$ and $d$.

The proof of the next proposition is different than that of Proposition \ref{transroot} (it also works for the latter). 

\begin{proposition}
\label{transfund}
$\mathbf{P\Gamma L}_{a + b}(k)$ acts transitively on the fundaments with ends of $\PG(a + b - 1,k) = \pi$ of type $(a,b;c,d)$. In particular, $\mathbf{P\Gamma L}_{a + b}(k)$ acts transitively on the fundaments of $\PG(a + b - 1,k)$ of type $(a,b)$.
\end{proposition}
{\em Proof}.\quad
Let $(\alpha,A,xy,\beta,B)$ and $(\alpha',A',x'y',\beta',B')$ be two fundaments with ends, both of type $(a,b;c,d)$, both in $\pi$. Let $(x,x_1,\ldots,x_{a - 1},y,y_1,\ldots,y_{b - 1})$ be an  ordered base of $\pi$ such that 
\begin{itemize}
\item
$(x,x_1,\ldots,x_{a - 1})$ is an ordered base of $\alpha$ and $(y,y_1,\ldots,y_{b - 1})$ an ordered base of $\beta$;
\item
$(x_{a - c},\ldots,x_{a - 1})$ is an ordered base of $A$ and $(y_{b - d},\ldots,y_{b - 1})$ is an ordered base of $B$.
\end{itemize}
Define in a similar way an ordered base $(x',x_1',\ldots,x_{a - 1}',y',y_1',\ldots,y_{b - 1}')$ with respect to $(\alpha',A',x'y',\beta',B')$.
Then $\mathbf{PGL}_{a + b}(k)$ contains an element sending the first ordered base to the second, as it acts transitively on the ordered bases of $\PG(a + b - 1,k)$.
\eop \\

\begin{corollary}
$\mathbf{PGL}_{a + b}(k)$ acts transitively on the fundaments of $\PG(a + b - 1,k)$ of type $(a,b)$. \eop \\
\end{corollary}

Let $\Gamma$ be the connected loose graph on two inner vertices $x$ and $y$, respectively of degree $b$ and $a$. We will show that 
\begin{equation}
\Aut(\mF(\Gamma) \times_{\Fun}k) \cong \Aut^{\mathrm{proj}}(\mF(\Gamma) \times_{\Fun}k)
\end{equation}
for any field $k$, where again $\Aut(\cdot)$ denotes the combinatorial group.

Suppose $c \leq a - 1$ edges on $y$  different from $xy$ have an end point, and that $d \leq b - 1$ edges on $x$ different than $xy$ have end points. 

Let $\A_x$ and $\A_y$ be the affine spaces
corresponding respectively to  the vertices $x$ and $y$ through the functor $\mF$. Write $X$ for $\mF(\Gamma)$, and $\mX_k$
for $\mF(\Gamma) \times_{\Fun} k$. 

As before, we want to see $\mX_k$ coming together with its embedding 
\begin{equation}
\mX_k \ \ \hookrightarrow\ \ \PG(a + b - 1,k).
\end{equation}

We projectively complete $\A_x$ and $\A_y$; the space of infinity of $\A_x$ is $\langle \beta_x,y \rangle$, with $\beta_x$ the projective space defined by the 
edges on $x$ different from $xy$, and the space at infinity of $\A_y$ is $\langle \alpha_y,x \rangle$, with $\alpha_y$ the projective space defined by the 
edges on $y$ different from $xy$. Also, let $A$ be the projective subspace of $\alpha_y$ defined by the end points ``in'' $\alpha_y$, and let $B$ be 
the projective subspace of $\beta_x$ defined by the end points in $\beta_x$. 
 
Any 
element of $\Aut^{\mathrm{proj}}(\mX_k)$ also fixes the projective completion $\overline{\mX_k}$. 

\begin{proposition}
We have that $\alpha$ is in $\Aut^{\mathrm{proj}}(\mX_k)$ (with $\alpha \in \mathbf{P\Gamma L}_{a + b}(k)$) if and only if 
$\alpha$ stabilizes the  incidence geometry of the fundament $(\alpha_y,A,xy,\beta_x,B)$.\eop 
\end{proposition}

The following is immediate.

\begin{proposition}
\label{kernelab}
The kernel of the action of $\mathbf{P\Gamma L}_{a + b}(k)_{\mX_k}$ on $\mX_k$ is trivial.
\end{proposition}

{\em Proof}.\quad
Let $\gamma \in \mathbf{P\Gamma L}_{a + b}(k)_{\mX_k}$ fix all the ($k$-rational) points of $\mX_k$. Then $\gamma$ fixes $\Pi_x := \overline{\A_x}$ and 
$\Pi_y := \overline{\A_y}$ pointwise. Consider any point $z$ in $\PG(a + b - 1,k)$ outside $\mX_k$. Then $\langle \Pi_x,z\rangle$ is a $(b + 1)$-dimensional projective space which meets the $a$-space $\Pi_y$ in a plane $\rho$ containing $xy$, and not contained in $\Pi_x$. So $\langle \Pi_x,z \rangle = \langle \Pi_x,\rho\rangle$. Hence
\begin{equation}
\langle \Pi_x,z \rangle^{\gamma} =  \langle \Pi_x,\rho\rangle^{\gamma} = \langle \Pi_x^{\gamma},\rho^{\gamma}\rangle = \langle \Pi_x,\rho\rangle, 
\end{equation} 
and the latter is pointwise fixed by $\gamma$, since $\Pi_x$ and $\rho$ are. (If $\gamma$ fixes $\mX_k$ pointwise, it also fixes each local affine space pointwise, so also their completions.)
So $z^{\gamma} = z$, and $\gamma$ is the identity.
\eop \\

For further purposes, let $\pi_x = \langle \beta_x,y\rangle$ the space at infinity of $\A_x$, and $\pi_y = \langle \alpha_y,x\rangle$ the space at infinity of $\A_y$.

The next couple of results carry over from roots to fundaments in a straightforward way.

\begin{proposition}
Let $E := \Aut(\Pi_x)_{(\pi_y,x,xy,y)}$ be the elementwise stabilizer of $\{\pi_y,x,y,xy\}$ in $\Aut(\Pi_x)$. (Here, $\Aut(\Pi_x)$ is the combinatorial automorphism group of $\Pi_x$, isomorphic to $\mathbf{P\Gamma L}_{a + 1}(k)$, and it is induced by $\mathbf{P\Gamma L}_{a + b}(k)$.)
Then each element of $E$ extends to an element of
$\PGaL_{a + b}(k)_{\mX_k}$ (in a not necessarily unique fashion).\eop
\end{proposition}

\begin{theorem}
Let $\Gamma$ be the loose graph defined in the beginning of this section. Then
\begin{equation}
\Aut(\mF(\Gamma) \times_{\Fun}k) \cong \Aut^{\mathrm{proj}}(\mF(\Gamma) \times_{\Fun}k)
\end{equation}
for any field $k$. \eop
\end{theorem}

If $S$ is a set of points in $\PG(a + b - 1,k)$, $\PGaL_{a + b}(k)_{[S]}$ denotes its pointwise stabilizer.

\begin{lemma}
Define $F := \Aut^{\mathrm{proj}}(\mX_k) \cap \PGaL_{a + b}(k)_{[\pi_y]}$, and  $G := \Aut^{\mathrm{proj}}(\mX_k) \cap \PGaL_{a + b}(k)_{[\pi_x]}$.
Then $F \cong \Aut(\Pi_y) \cap \PGL_{a + b}(k)$, and   $G \cong \Aut(\Pi_x) \cap \PGL_{a + b}(k)$. \eop
\end{lemma}

The following theorem is proved in exactly the same way as in the case of roots.

\begin{theorem}
Let $\PGL(\mX_k)_{(x,y)}$ be defined as 
\begin{equation}
\Aut^{\mathrm{proj}}(\mX_k)_{(x,y)} \cap \PGL_{a + b}(k).
\end{equation}
Then $\PGL(\mX_k)_{(x,y)}$ is isomorphic to the internal central product of $F$ and $G$. \eop \\
\end{theorem}

The general version is the following.

\begin{theorem}[Trees on two inner vertices]
We have that
\begin{equation}
\Aut^{\mathrm{proj}}(\mX_k)\ \cong\ \Aut(\mX_k)\ \cong\ ((F * G) \rtimes \Aut(k)) \rtimes \langle \varphi \rangle. 
\end{equation}
In the latter expression, $\varphi$ is trivial, unless the type of the fundament has the form $(a,a;c,c)$, in which case $\varphi$ is 
an involution in the automorphism group of $\Gamma$ which switches $x$ and $y$. 
\end{theorem}
{\em Proof}.\quad
If the type is $(a,a;c,c)$, then obviously there is an involution as in the statement of the theorem. And any element in $\Aut(\Gamma)$
fixes both $x$ and $y$ if $a \ne b$ or $c \ne d$. \eop

\medskip
\subsection{General loose trees}

Let $T = (V,E)$ be a finite loose tree, and assume its number of vertices is at least $3$. Let $\overline{T}$ be the graph theoretical completion of $T$; define the {\em boundary} of $T$, denoted $\partial(T)$, as the set of vertices of degree $1$ in $\overline{T}$. Let $x$ be a vertex which is at distance $1$ from $\partial(T)$ (i.e., is adjacent with at least one vertex of $\partial(T)$). As $\vert V \vert \geq 3$, $x$ is an inner vertex of degree at least $2$.  

Define $k$ and $\mX_k$ as before. Let $\mathbf{PG}(m - 1,k)$ be the ambient projective space of $\mX_k$.

\begin{proposition}
\label{propkernel}
The kernel of the action of $\mathbf{P\Gamma L}_{m}(k)_{\mX_k}$ on $\mX_k$ is trivial.
\end{proposition}

{\em Proof}.\quad
Let $\gamma \in \mathbf{P\Gamma L}_{m}(k)_{\mX_k}$ fix all the $k$-rational points of $\mX_k$. If $T$ is an affine $\Fun$-space (with some end points), then
there is nothing to prove. So suppose $T$ is not.

Define $\Pi_x := \overline{\A_x}$ as before, and let $y \sim x \ne y$ be not in $\partial(T)$ (such a point exists). Let $\Pi_y$ be the projective completion of the loose graph $T_y$
induced on the vertex set $V_y := \{ v \in V \vert d(v,x) \geq 2 \} \cup \{ y \}$ (by ``induced,'' we mean, besides inheriting the induced {\em graph} structure, that if $e$ is a loose edge in $T$ which is incident with 
a vertex of $V_y$, then $e$ is in $T_y$). Note that $\Pi_y$ contains $\overline{\A_y}$.

Now repeat the argument of Proposition \ref{kernelab}, using induction on the loose tree $T_y$. \eop \\

In the next couple of results, we keep using the notation introduced in the beginning of this subsection. Also, with $I$ the set of inner
vertices of $\overline{T}$, and $w \in I$, let $S(w)$ be the subgroup of $\Aut^{\mathrm{proj}}(\mX_k)$ which fixes the $k$-rational points of $\mX_k$ inside all affine subspaces $\widetilde{\A_v}$ (over $k$) which are generated (over $\Fun$) by a vertex $v$ different from $w$ and all directions on $v$ which are not incident with $w$.  So, if the distance of $v$ to $w$ is at least $2$, the local space at $v$ is fixed pointwise, and if the distance is $1$, $\widetilde{\A_v}$ is an affine space of dimension one less than the dimension of $\A_v$.
(In particular, the points in $I \cap \mathbf{B}(w,1)$ are fixed.) 


In the next theorem, one recalls that $\mX_k$ comes with an embedding
\begin{equation}
T \ \hookrightarrow \ \mX_k \ \hookrightarrow\ \PG(m - 1,k),
\end{equation}
so that it makes sense to consider stabilizers of substructures of $T$ in, e.g.,  $\PGL(\mX_k)$.\\

\begin{theorem}
\label{cenprod}
Let $\PGL(\mX_k)_{[I]}$ be defined as 
\begin{equation}
\Aut^{\mathrm{proj}}(\mX_k)_{[I]}\ \cap\ \PGL_{m}(k).
\end{equation}
Then $\PGL(\mX_k)_{[I]}$ is isomorphic to 
\begin{equation}
\prod^{\mathrm{centr}}_{w \in I}S(w).
\end{equation}
\end{theorem}
{\em Proof}.\quad
Let $x \in I$ be at distance $1$ from $\partial(T)$. Also, let $y \sim x \ne y$, $y \not\in \partial(T)$ and $y \in I$. Let $T_y$ be the loose graph
induced on the vertex set $V_y := \{ v \in V \vert d(v,x) \geq 2 \} \cup \{ y \}$, but without the edge $xy$. Let $H(y)$ be the subgroup of $\Aut^{\mathrm{proj}}(\mX_k)$ which fixes
the affine subspace of $\PG(m - 1,k)$ pointwise that is generated by all edges on $x$ in $T$ except $xy$. It is important to observe that $S(y) \leq H(y)$ for the induction argument later on.
Then in the same way as in the proof of Theorem \ref{thmcp},
one shows that 
\begin{equation}
\PGL(\mX_k)_{[V]} = S(x) \ast H(y).
\end{equation}
Now perform induction on $T_y$ to conclude that 
\begin{equation}
\PGL(\mX_k)_{[I]} = S(x) \ast \Big( S(y) \ast \Big( \cdots \Big)\Big).
\end{equation}
Note that for each $u, v \in I$, it follows that
\begin{equation}
[S(u),S(v)] \ = \ \{\id\}.
\end{equation}
\eop

\medskip
\subsubsection{Determination of $S(w)$}

We start by remarking that although in general $S(w)$ fixes a lot of points, it is not necessarily a subgroup of $\PGL_m(k)$ (see for instance Lemma 
\ref{lemfield} below). What we {\em do} know | by its mere definition | is that it is a subgroup of $\PGaL_{m}(k)$.

We will distinguish two cases in order to determine $S(w)$.

\subsection*{\quad $\dagger$ $w$ is the only inner point}

Then all the edges are incident with $w$. 
Call $E$ the set of such edges with an end point, and $L$ the set of loose edges. Put $\vert E \vert = e$ and $\vert L \vert = \ell$. Then obviously
\begin{equation}
S(w) \cong {\Big(\PGaL_{e + \ell + 1}(k)_{L}\Big)}_{[E \cup \{w\}]}.
\end{equation}
By the first remark of this subsection, it is not contained in the projective linear subgroup.

\subsection*{\quad $\ddagger$ $w$ is not the only inner point}

Then there is some inner vertex $v \sim w$ different from $w$ which is itself incident  to some edge $W \ne wv$. 
Now over $k$, the projective line which is the completion of the affine line determined by the incident vertex-edge pair $(v,W)$, is fixed pointwise 
by $S(w)$, so $S(w)$ must be a subgroup of $\PGL_m(k)$. 

Let $E$ be the set of edges incident with $w$ which have an end point, let $L$ be the set of loose edges incident with $w$, and let $I$ be the 
set of edges incident with $w$ which are incident with another inner point. Put $\vert E \vert = e$, $\vert L \vert = \ell$ and $\vert I \vert = i$. Let $\delta$ 
be an element of $S(w)$; then it induces an element of $\PGL(\overline{\A_w})$ (the latter meaning the projective linear group of the local 
projective space at $w$). If $\delta'$ is another such element which induces the same action, it is obvious that $\delta\delta^{-1}$ is the identity on 
the entire ambient space $\PG(m - 1,k)$. So $S(w)$ faithfully is a subgroup of ${\Big(\PGL_{e + \ell + i + 1}(k)_L\Big)}_{[E \cup I \cup \{w\}]}$.

Note that the projective space generated (over $k$) by the points at distance at least $2$ from $w$ in $\overline{\Gamma}$, is fixed pointwise by $S(w)$.
So in particular $\pi_I$, the projective space generated by the inner vertices adjacent to $w$, is also fixed pointwise. It now follows easily that 
\begin{equation}
S(w) \cong {\Big(\PGL_{e + \ell + i + 1}(k)_L\Big)}_{[E \cup \pi_I \cup \{w\}]}.
\end{equation}

\medskip
\subsubsection{Caution: central and direct products}

On the graph theoretical level (that is, on the combinatorial $\Fun$-level), the groups which occur in Theorem \ref{cenprod} are much easier to describe, 
replacing the central product by a direct product. The central product is needed as soon as $k^\times$ is not trivial.

\medskip
\subsubsection{Inner Tree Theorem}

The following theorem is a crucial ingredient in the proof of our main theorem for trees. 

\begin{theorem}[Inner Tree Theorem]
\label{innertree}
Let $T$ be a loose tree, and let $k$ be any field. As usual, put $\mX_k = \mF(T) \otimes_{\Fun} k$, and consider the embedding
\begin{equation}
\iota: T \ \hookrightarrow\ \mX_k.
\end{equation}
Let $\Aut(\mX_k)$ be any of the automorphism groups which are considered in this paper (i.e., combinatorial, induced by projective space, topological, or scheme-theoretic). 
Let $I$ be the set of inner vertices of $T$, and let $T(I)$ be the subtree of $T$ induced on $I$. 
Then if $\vert I \vert \geq 2$, we have that $\Aut(\mX_k)$ stabilizes $\iota(T(I))$. Moreover, $\Aut(\iota(T(I)))$ is induced by $\Aut(\mX_k)$.
\end{theorem}

{\em Proof}.\quad
Each edge of $\iota(T(I))$ defines a projective line over $k$ which is a full line of the ambient space of $\mX_k$. Let $\iota(T(I))_k$ be this set of projective lines.
Now define $\underline{\mX_k}$ as the {\em projective part} of $\mX_k$ | by definition, it is the union of all projective $k$-lines which are completely contained in $\mX_k$. As each local affine space at a vertex of $T$ is an affine space with some possible end points at infinity, one observes that $\underline{\mX_k}$ consists precisely of the projective $k$-lines which are defined by the edges with two different vertices of $T$. That is, $\underline{\mX_k}$ consists of $\iota(T(I))_k$ together with 
additional projective lines defined by edges which contain both an inner vertex and an end point of $T$. 
As $\vert I \vert \geq 2$, the first part of the theorem easily follows. 

That $\Aut(\iota(T(I)))$ is induced follows by functoriality (and the discussion in \S \ref{diffaut}).
\eop \\

Note that if $\vert I \vert = 1$, $T$ defines an affine space with some end points, so the theorem is not true, unless its dimension is $0$. If $\vert I \vert = 0$, then either $T$ is the empty tree, or $T$ is an edge with one or two vertices.

\medskip
\subsubsection{The general group}

Before proceeding, we need another lemma. We use the notation of the previous subsection.

\begin{lemma}[Field automorphisms]
\label{lemfield}
Let $\PG(m - 1,k)$ be the ambient space of $\mX_k$. 
We have that 
\begin{equation}
{\PGaL_{m}(k)}_{\mX_k}\Big/{\PGL_m(k)}_{\mX_k} \cong k^{\times}.
\end{equation}
\end{lemma}
{\em Proof}.\quad
Let $\Delta$ be the base of $\PG(m - 1,k)$ corresponding to the vertices of $\overline{T}$. Then it is well known that 
\begin{equation}
{\PGaL_{m}(k)}_{[\Delta]}\Big/{\PGL_m(k)}_{[\Delta]} \cong k^{\times}.
\end{equation}
(In fact, working with homogeneous coordinates with respect to $\Delta$, ${\PGaL_m(k)}_{[\Delta]}$ contains all elements of the form
$\mathbf{x} \mapsto \id_m\mathbf{x}^\tau$, with $\mathbf{x}$ a column vector representing points in homogeneous coordinates, $\id_m$ the identity $(m \times m)$-matrix
and $\tau \in \Aut(k)$.) As ${\PGaL_{m}(k)}_{[\Delta]} \leq {\PGaL_{m}(k)}_{\mX_k}$ and ${\PGL_{m}(k)}_{[\Delta]} \leq {\PGL_{m}(k)}_{\mX_k}$, the lemma easily follows. \eop \\

\begin{theorem}[Projective automorphism group]
\label{MTtrees}
Let $T$ be a loose tree, and let $k$ be any field. Put $\mX_k = \mF(T) \otimes_{\Fun} k$, and consider the embedding
\begin{equation}
\iota: T \ \hookrightarrow\ \mX_k.
\end{equation}
Let $I$ be the set of inner vertices of $T$, and let $T(I)$ be the subtree of $T$ induced on $I$. 
We have $\PGaL(\mX_k) = \Aut^{\mathrm{proj}}(\mX_k)$ is isomorphic to 
\begin{equation}
\Big(\Big(\prod^{\mathrm{centr}}_{w \in I}S(w)\Big) \rtimes \Aut(T(I))\Big) \rtimes k^{\times}.
\end{equation}
\end{theorem}
{\em Proof}.\quad
First note that by Proposition \ref{propkernel}, the kernel of the action of ${\PGaL_m(k)}_{\mX_k}$ on $\mX_k$ is trivial.
Then by Lemma \ref{lemfield}, we only have to show that 
\begin{equation}
{\PGL_m(k)}_{\mX_k}\ \cong \ \Big(\prod^{\mathrm{centr}}_{w \in I}S(w)\Big) \rtimes \Aut(T(I)).
\end{equation}
By Theorem \ref{cenprod}, we have that $\PGL(\mX_k)_{[I]}$ is isomorphic to 
\begin{equation}
\prod^{\mathrm{centr}}_{w \in I}S(w),
\end{equation}
and obviously $\PGL(\mX_k)_{[I]} \unlhd {\PGL_m(k)}_{\mX_k}$. 

The theorem now follows from the Inner Tree Theorem. 
\eop \\

\subsubsection{The combinatorial automorphism group}

By Theorem \ref{MTtrees}, we can now determine the combinatorial group as well. 

\begin{theorem}[Combinatorial automorphism group]
\label{autcombloose}
Let $T$ be a loose tree, and let $k$ be any field. Put $\mX_k = \mF(T) \otimes_{\Fun} k$, let $I$ be the set of inner vertices, and suppose that $\vert I \vert \geq 2$. Let $\iota$ be as in Theorem \ref{MTtrees}.
Then
\begin{equation}
\Aut^{\mathrm{comb}}(\mX_k) \cong \Aut^{\mathrm{proj}}(\mX_k).
\end{equation}
\end{theorem}
{\em Proof}.\quad
As in the proof of Theorem \ref{6.7}, we assume by way of contradiction that there is some $\alpha \in \Aut^{\mathrm{comb}}(\mX_k) \setminus \Aut^{\mathrm{proj}}(\mX_k)$. As  in that theorem, by the fact that $\Aut^{\mathrm{proj}}(\mX_k)$ induces $\Aut(\iota(T(I)))$ by the Inner Tree Theorem, we may assume that $\alpha$ fixes all vertices of $\iota(T(I))$. Now  $\alpha$ induces projective automorphisms in each $\overline{\A_x}$ with $x$ an inner vertex, which are compatible on edges of $\iota(T(I))$. By Theorem \ref{MTtrees}, we can end in the same way as in the proof of Theorem \ref{6.7}. 
\eop \\

\medskip
\section{More on the different automorphism group types}
\label{diffaut}

We have shown in Proposition \ref{combtopo} that for each $\mX_k$, the combinatorial automorphism group is a subgroup of the topological automorphism group. Also it is clear that any projectively induced automorphism is combinatorial, but the other direction is in general {\em not} true. 
Let $\Gamma$ be, e.g., an edge with two different vertices, so that for all $k$, $\mX_k$ is a projective $k$-line. Then each permutation of the $k$-points yields a combinatorial automorphism, but not all of these come from projective automorphisms for all $k$. 
So
\begin{equation}
\begin{cases}
\Aut^{\mathrm{top}}(\mX_k) \geq \Aut^{\mathrm{comb}}(\mX_k)\\
\Aut^{\mathrm{comb}}(\mX_k), \Aut^{\mathrm{top}}(\mX_k) \geq \Aut^{\mathrm{proj}}(\mX_k).
\end{cases}
\end{equation}

Any projective automorphism is also scheme-theoretic, but not the other way around | think of affine spaces as a typical example. Any scheme-theoretic automorphism also induces a topological one, but in some cases this might not happen in a faithful way, and this is an aspect we want to come back to in a subsequent paper. Not every scheme-theoretic automorphism is combinatorial (think again of affine spaces), but the {\em linear} ones are.
In any case,
\begin{equation}
\begin{cases}
\Aut^{\mathrm{sch}}(\mX_k)\Big/ (\text{possible kernel}) \leq \Aut^{\mathrm{top}}(\mX_k)\\
\Aut^{\mathrm{sch, linear}}(\mX_k)\Big/ (\text{possible kernel}) \leq \Aut^{\mathrm{comb}}(\mX_k).
\end{cases}
\end{equation}

We do not know yet whether every combinatorial automorphism (say, in sufficiently high dimension) is of scheme-theoretic origin.

\medskip
\section{Convexity}

Let $T$ be a loose tree, and for any field $k$, consider $\mX_k := \mF(T) \otimes_{\Fun} k$. In this section we will prove a useful convexity property
for the spaces $\mX_k$.

The following lemma is trivial, but also useful.

\begin{lemma}
\label{lemsubgraph}
Let $G$ be any subgraph of $\overline{T}$, not necessarily connected. Then the dimension of the projective space generated over $\Fun$ by $G$ equals
the number of vertices of $G$ minus $1$. \eop
\end{lemma}

\begin{theorem}[Convexity]
Let $k$ and $\mX_k$ be as in the beginning of this section. Let $\A_u$ and $\A_v$ be local affine spaces over $k$  with $u, v \ne u$ vertices of $T$. If $x \in \A_u$, but 
not contained in any of the lines determined by the local loose star of $u$, and $y \in \A_v$ is not contained in any of the lines determined by the local loose star of $v$,
then the projective $k$-line $xy$ only meets $\mX_k$ in $x$ and $y$. 
\end{theorem}
{\em Proof}.\quad
Suppose by way of contradiction that $z \in (\mX_k \cap xy) \setminus \{x,y\}$. Obviously $z \not\in \A_u \cup \A_v$, so $z$ is in some other local affine $k$-space $\A_w$, with $w$ a vertex of $T$. 
There are (essentially) five possible configurations to be considered:
\begin{itemize}
\item[(1)]
$u \sim v$ and $u \not\sim w \not\sim v$;
\item[(2)]
$u \sim v$ and $u \not\sim w \sim v$;
\item[(3)]
$u \not\sim v$ and $u \not\sim w \not\sim v$;
\item[(4)]
$u \not\sim v$ and $u \not\sim w \sim v$;
\item[(5)]
$u \not\sim v$ and $u \sim w \sim v$.
\end{itemize}
Note that $u, v, w$ cannot form a triangle. In each of the cases, consider the projective space generated by $\A_u, \A_v$ and $\A_w$, calculate its dimension, and and apply Lemma \ref{lemsubgraph} to find a contradiction.
\eop

\bigskip
\section{The edge-relation dichotomy}

The fact that the calculations for loose trees $T$ are so successful rests largely on the fact that there are no cycles; that property leads to the fact that we can apply the Inner Tree Theorem, and this makes it able to determine the various automorphism groups of $\mF(T) \otimes k$, $k$ any field.

The examples which are the farthest from satisfying the Inner Tree Theorem are affine and projective spaces. In case of affine spaces $\A_k^n$, the automorphism group (be it combinatorial or scheme-theoretic) acts transitively on the $k$-points, so obviously the Inner Tree Theorem, formulated for loose graphs (see \S \ref{InnGrProp}),  cannot hold. In fact, we have the following observation the trivial proof of which we leave to the reader.

\begin{theorem}
Let $\Gamma$ be the loose graph of an affine or projective $\Fun$-space. Then for any field $k$ and any of the considered automorphism groups $\Aut(\cdot)$, we have
that $\Aut(\mF(\Gamma)\otimes k)$ acts transitively on the set of subgeometries isomorphic to $\Gamma$. (Here, as before a subgeometry consists of $k$-points and affine or projective $k$-lines.) \eop 
\end{theorem}

\medskip
\subsection{Examples close to trees}

Consider the following loose graph $\Gamma_1$, which, for each field $k$, defines in the ambient projective $3$-space $\PG(3,k)$,  four affine planes each with two extra points at infinity and cyclically denoted by $\alpha_i$ ($i = 1,2,3,4$), in which ``adjacent planes'' meet in a projective line, and ``opposite planes" meet precisely in the end points. 
Denote the scheme by $\mX_k$.

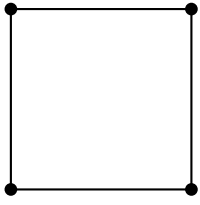
\begin{figure}[h]\label{p1}

\begin{tikzpicture}[style=thick, scale=1.2]

\draw (0,0)-- (2,0);

\draw (0,0)-- (0,2);

\draw (0,2)-- (2,2);

\draw (2,2)-- (2,0);

\fill (0,0) circle (2pt);

\fill (2,2) circle (2pt);

\fill (0,2) circle (2pt);

\fill (2,0) circle (2pt);

\end{tikzpicture}

\caption{The loose graph $\Gamma_1$}

\end{figure}

Obviously we have

\begin{equation}
\Aut^{\text{proj}}(\mX_k)\ \cong\ {\PGaL_4(k)}_{\Gamma_1},
\end{equation}
where $\Gamma_1$ comes with the embedding

\begin{equation}
\Gamma_1\ \hookrightarrow\ \mX_k.
\end{equation}

The complement $\Gamma_1^c$ of $\Gamma_1$ in its ambient projective $\Fun$-space is also fixed by $\Aut^{\text{proj}}(\mX_k)$, as that complement just defines two disjoint multiplicative groups. Notice however that 

\begin{equation}
\Big(\mF(\Gamma_1) \otimes_{\Fun} k\Big) \ \coprod \ \Big(\mF(\Gamma_1^c) \otimes_{\Fun} k\Big)\ \ne \ \PG(3,k)!
\end{equation}

The example $\Gamma_1$ easily generalizes to the class of polygonal graphs $\Gamma(m)$ with $m + 1$ vertices, $m \geq 0$, $m \ne 2$; for $m = 0,1$ we get 
spaces $\Proj(k)$ and $\Proj(k[X])$ which satisfy the Inner Tree Property; for $m = 3$ we get a projective $k$-plane, and for $m \geq 4$, we obtain a scheme consisting of $m + 1$ affine $k$-planes each with two extra points at infinity, which intersect two by two according to their graph intersection (in a projective $k$-line, a point or no intersection). All of them except $\Gamma(2)$  have the property that
\begin{equation}
\Aut^{\text{proj}}(\mX_k)\ \cong\ {\PGaL_4(k)}_{\Gamma(m)}.
\end{equation}

The graph complements are also fixed by $\Aut^{\text{proj}}(\mX_k)$.

\medskip
\subsection{Missing piece}

Let $\Gamma$ be a loose graph, $k$ any field, $\hP_k := \PG(m - 1,k)$ the ambient space over $k$, and $\Gamma^c$ the complement in $\hP_{\Fun}$ of 
$\Gamma$. We have a decomposition

\begin{equation}
\Big(\mF(\Gamma) \otimes_{\Fun} k\Big) \ \coprod \ \Big(\mF(\Gamma^c) \otimes_{\Fun} k\Big)\  \ \coprod \ \mY_k(\Gamma) \ = \ \PG(m - 1,k),
\end{equation}
for some (quasi-projective) variety $\mY_k(\Gamma)$. The variety $\mY_k$ measures a difference in behavior of $\mF(\Gamma) \otimes_{\Fun} k$ with respect to fields $k$ and $k = \Fun$, since, for instance, for $k = \Fun$ we have that $\mF(\Gamma) \ \coprod \ \mF(\Gamma^c)$ partitions the line set of $\PG(m - 1,\Fun)$. (Note however that one has to be careful with decompositions in terms of loose graphs: e.g., an affine $\Fun$-plane minus a multiplicative group $\bG_m^1$ is {\em not} an affine line! | one might want to think in terms of the Grothendieck ring of $\Fun$-schemes $K_0(\Sch_{\Fun})$ \cite{MMKT} to see this more clearly.)

It might be interesting to study the maps
\begin{equation}
\mY_k: \Gamma \ \mapsto\ \mY_k(\Gamma).
\end{equation}

\medskip
\subsection{Examples close to the ambient space}

Now consider the following example $\Gamma_2$, which, for each field $k$, defines a projective $3$-space $\PG(3,k)$ without one multiplicative group $\bG_m$ (corresponding to the missing diagonal edge). (Denote the scheme by $\mX_k$.)

\begin{figure}[h]\label{p1}

\begin{tikzpicture}[style=thick, scale=1.2]

\draw (0,0)-- (2,0);

\draw (0,0)-- (0,2);

\draw (0,2)-- (2,2);

\draw (2,2)-- (2,0);

\draw (0,2)-- (2,0);

\fill (0,0) circle (2pt);

\fill (2,2) circle (2pt);

\fill (0,2) circle (2pt);

\fill (2,0) circle (2pt);

\end{tikzpicture}

\caption{The loose graph $\Gamma_2$}

\end{figure}
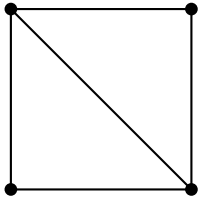

Let $x$ and $y$ be the two $k$-points of $\PG(3,k)$ in the projective line defined by $\bG_m$ which are not contained in $\bG_m$. Then obviously

\begin{equation}
\Aut^{\text{proj}}(\mX_k)\ \cong\ {\PGaL_4(k)}_{\{x,y\}},
\end{equation}
so $\Aut^{\text{proj}}\mX_k$ does not fix the graph defined by

\begin{equation}
\Gamma_2\ \hookrightarrow\ \mX_k.
\end{equation}

What it {\em does} fix, is the complement of $\Gamma_2$ in the projective $\Fun$-space defined by $\Gamma_1$ (considered in the same embedding).

\medskip
\subsection{Schemes satisfying the Inner Graph Property}
\label{InnGrProp}

One essential ingredient in the proof of our main theorem for trees, is the {\em inner tree property}, which we define as follows for general loose graphs.

\begin{quote}
Let $\Gamma$ be a loose graph, and let $k$ be any field. Put $\mX_k = \mF(\Gamma) \otimes_{\Fun} k$, and consider the embedding
\begin{equation}
\iota: \Gamma \ \hookrightarrow\ \mX_k.
\end{equation}
Let $\Aut(\mX_k)$ be one of the automorphism groups considered in this paper | combinatorial, induced by projective space or scheme- theoretic. 
Let $I$ be the set of inner vertices of $T$, and let $\Gamma(I)$ be the subgraph of $\Gamma$ induced on $I$. Suppose $\vert I \vert \geq 2$. Then we say that $\Gamma$ satisfies the {\em inner graph property} if $\Aut(\mX_k)$ stabilizes $\iota(\Gamma(I))$. \\
\end{quote}

\begin{question}
Characterize (the) loose graphs that do/do not have the inner graph property.
\end{question}

Let $\wis{InnGraph}$ be the category of loose graphs
which have the inner graph property. Following the same lines of the proof of Theorem \ref{MTtrees}, one can determine the map

\begin{equation}
\Aut: \wis{InnGraph}\ \mapsto\ \wis{Group}:\ \Gamma\ \mapsto\ \Aut(\Gamma).
\end{equation}

We will handle this case in a forthcoming paper \cite{MMKTfunct2}.\\

\medskip
\subsection{Heisenberg principle}

Let $\wis{LGraph}$ be the category of loose graphs, $\wis{LTree}$ the category of loose trees, and $\wis{CGraph}$ the category of complete graphs. 
We end the paper with the following questions.

\begin{question}
Does there exist a distance function 
\begin{equation}
\delta: \wis{LGraph} \times {\wis{LGraph}}\ \mapsto\ (S,\leq),
\end{equation}
with $(S,\leq)$ a (totally) ordered set, such that the following properties hold?

\begin{itemize}
\item
The distance between a loose tree and its completion in $\wis{CGraph}$ is maximal.
\item
If $\mathrm{min}\{\delta(\Gamma,T)\ \vert\ T \in \wis{LTree}, T \leq \Gamma \}\ \ll$, then $\Gamma$ satisfies the inner graph property.
\item
If $\delta(\Gamma,\overline{\Gamma})\ \ll$, with $\overline{\Gamma}$ the completion of $\Gamma$ in $\wis{CGraph}$, then $\Gamma$ does not satisfy the inner graph property.\\
\end{itemize}
\end{question}

We strongly suspect that $\delta$ should be expressed in terms of cycles. 

\begin{question}
Let $\delta$ be as in the previous question. Let $\Gamma$ be in $\wis{LGraph}$, and suppose that
\begin{equation}
\mathrm{min}\{\delta(\Gamma,T)\ \vert\ T \in \wis{LTree}, T \leq \Gamma \} \cdot \delta(\Gamma,\overline{\Gamma})
\end{equation}
is ``quadratic,'' when can one decide that $\Gamma$ satisfies the inner graph property? 
\end{question}

\end{document}